\documentclass[9pt,twocolumn,twoside]{pnas-new}
\usepackage{geometry}                
\usepackage{graphicx}
\usepackage{epstopdf}
\templatetype{pnasresearcharticle} 
\usepackage{booktabs,amsfonts,amsmath,mathtools,framed,pbox,graphicx,caption,amssymb,subcaption}
\usepackage{xcolor}

\newcommand{\re}{\color{black} }
\newcommand{\ed}{\color{black} }

\usepackage{caption}
\usepackage[labelformat=simple]{subcaption}

\newcommand{\ie}{\textit{i.e. }}
\newcommand{\eg}{\textit{e.g. }}

\title{Diversity-enabled sweet spots in layered architectures and speed-accuracy trade-offs in sensorimotor control} 

\author
{Yorie Nakahira$^{a,e}$, Quanying Liu$^{b,e}$, Terrence J. Sejnowski$^{c,d,\ast}$, John C. Doyle$^{e,\ast}$\\
\normalsize{$^{a}$College of Engineering, Carnegie Mellon University, Pittsburgh, PA, USA}\\ 
\normalsize{$^{b}$Department of Biomedical Engineering, Southern University of Science and Technology, Shenzhen, Guangdong, China }\\
\normalsize{$^{c}$The Salk Institute for Biological Studies, La Jolla, CA, USA}\\
\normalsize{$^{d}$Division of Biological Sciences, University of California, San Diego, La Jolla, CA, USA}\\
\normalsize{$^{e}$Division of Engineering and Applied Science, California Institute of Technology,} \normalsize{Pasadena, CA, USA}\\
\normalsize{$^\ast$To whom correspondence should be addressed; E-mail: terry@salk.edu, doyle@caltech.edu.}
}

\significancestatement{Nervous systems use highly effective layered architectures in the sensorimotor control system to minimize the harmful effects of delay and inaccuracy in biological component. To study what makes effective architectures, we develop a theoretical framework that connects the component speed-accuracy trade-offs (SATs) with system SATs and characterizes the system performance of a layered control system. We show that diversity in layers (\eg planning and reflex) allows fast and accurate sensorimotor control even when each layer uses slow or inaccurate components. We term such phenomena ``diversity-enabled sweet spots (DESSs)''. DESSs explain and link the extreme heterogeneities in axon sizes and numbers at the component level and the resulting robust performance in sensorimotor control.
}

\authorcontributions{Please provide details of author contributions here.}
\authordeclaration{Authors declare no conflict of interest.}

\begin{abstract}

Nervous systems sense, communicate, compute and actuate movement using distributed components with severe trade-offs in speed, accuracy, sparsity, noise and saturation. Nevertheless, brains achieve remarkably fast, accurate, and robust control performance due to a highly effective layered control architecture. Here we introduce a driving task to study how a mountain biker mitigates the immediate disturbance of trail bumps and responds to changes in trail direction. We manipulated the time delays and accuracy of the control input from the wheel as a surrogate for manipulating the characteristics of neurons in the control loop. The observed speed-accuracy trade-offs (SATs) motivated a theoretical framework consisting of layers of control loops with components having diverse speeds and accuracies within each physical level, such as nerve bundles containing axons with a wide range of sizes. Our model explains why the errors from two control loops -- one fast but inaccurate reflexive layer that corrects for bumps, and a planning layer that is slow but accurate -- are additive, and show how the errors in each control loop can be decomposed into the errors caused by the limited speeds and accuracies of the components. 
These results demonstrate that an appropriate diversity in the properties of neurons across layers helps to create ``diversity-enabled sweet spots'' (DESSs) so that \textit{both} fast \textit{and} accurate control is achieved using slow or inaccurate components.

\end{abstract}

\dates{This manuscript was compiled on \today}
\doi{\url{www.pnas.org/cgi/doi/10.1073/pnas.XXXXXXXXXX}}

\begin{document}

\maketitle

\thispagestyle{firststyle}

\ifthenelse{\boolean{shortarticle}}{\ifthenelse{\boolean{singlecolumn}}{\abscontentformatted}{\abscontent}}{}

When riding a mountain bike down a twisting and bumpy trail, skilled riders can descend safely without crashing despite limitations imposed by imperfect components in the brain and trade-offs between traveling fast and staying on the trail. What enables such remarkably robust performance in complex and uncertain environments? 
Although this question is of great importance in both science and engineering, it has received little attention in the prior work in neuroscience and control. 

The remarkable robustness of sensorimotor control has fostered the widespread illusion that system performance is unconstrained by the limitations of its components~\cite{doyle2011architecture}. Consequently, little attention has been paid to understanding the design principles that deconstrain the limitations of its components. However, in both biological and engineered systems, ignoring the hard limits results in fragility and may even lead to catastrophic failure.

\begin{figure}[ht]
\center
\includegraphics[width=\linewidth]{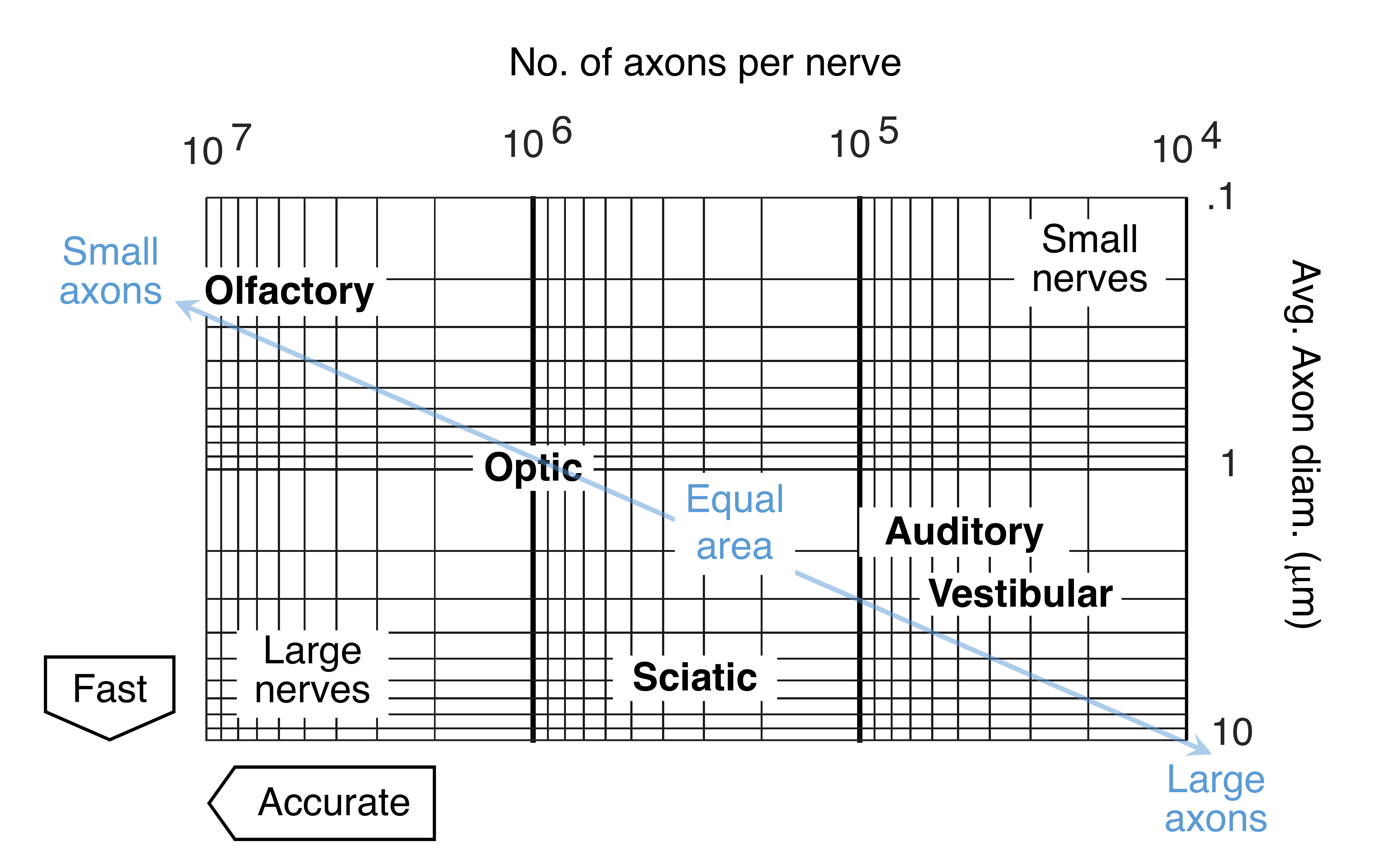}
\caption{\textbf{Sizes and numbers of axons for selected nerves and the resulting SATs.} The dashed line represents nerves with an equal cross-sectional area, which is proportional to $\lambda$ in Eq.~\ref{eq:delay-rate-spike1}--\ref{eq:spike-rate-code}. The nerves shown have similar cross-sectional areas but wildly different compositions of axon size and number, resulting in different speed and accuracy in nerve signaling~\cite{laughlin2003communication,sterling2015principles}.  
A myelin sheath around an axon can also increase its speed of propagation. 
Many nerves, such as the sciatic nerve, contain a mixture of axons with different sizes and degrees of myelination.}
 \label{fig:axons}
\end{figure}

A clue to this puzzle lies in the striking contrast of speed-accuracy tradeoffs (SATs) at the component level and SATs at the system level. The constraints on sensory and motor nerves that implement sensorimotor control are often stringent. For example, the spatial and metabolic costs to build, operate and maintain signaling in nerves constrain the fiber sizes and numbers of axons in a nerve. This limits the speed and the amount of information that these axons can transmit~\cite{laughlin2003communication,sterling2015principles}. Large nerves with axons that are both large in size and number are rare (Fig.~\ref{fig:axons}), which suggests that achieving both speed and accuracy may be prohibitively expensive. 

Such component limits could constrain the sensorimotor control to be slow and/or inaccurate in a naive design. However, in practice our cognitive decision-making and sensorimotor control are remarkably robust, fast, and accurate as if the component limits were deconstrained~\cite{standage2015toward,doya2007bayesian,chittka2009speed,card2018psychology,heitz2013neural,saxena2020performance}. Examples can be observed in the extraordinary performance of athletes, mountain biking among others, and power-laws in reaching such as Fitts' Law~\cite{card2018psychology,Nakahira2020Fitts}.

This striking contrast between system and component SATs in sensorimotor control suggests there are highly efficient mechanisms that successfully deconstrain component limitations in the sensorimotor system, so that component constraints are not apparent. {\re Strictly speaking, it is not possible for the aggregate information rate of all components to exceed the sum of the information rates of each component. Although the component constraints cannot be deconstrained at the component level, these constraints are unconstrained at the system level.} This can be achieved using virtualization and layered architectures in network engineering, and we show here that similar principles are found in brains as well.

There are two major challenges in understanding the design principles found in nature that deconstrain the component limits in sensorimotor control. The first is to bridge the SATs at the level of neurophysiology and the SATs at the level of system and behavior. The second is to understand the integration and coordination of layers with distinct roles with heterogeneous components and limitations. Despite extensive research focused on individual levels and layers, there are few theoretical insights or experimental tools available to integrate the component constraints of individual layers with fundamental limits on the performance of the entire system. On the theory side, we do not yet know enough about neural coding and control mechanisms to establish a complete model for control pathways from sensory to motor units or to pinpoint performance bottlenecks. On the experimental side, it is difficult to noninvasively manipulate the properties of the components, including time delays and information rates, to observe how they influence the system SATs.  

\begin{figure}[hb!]

\begin{subfigure}[b]{\linewidth}
        \caption{Sensorimotor control system used in oculomotor control for visual tracking}
         \centering
        	\includegraphics[width=0.8\linewidth]{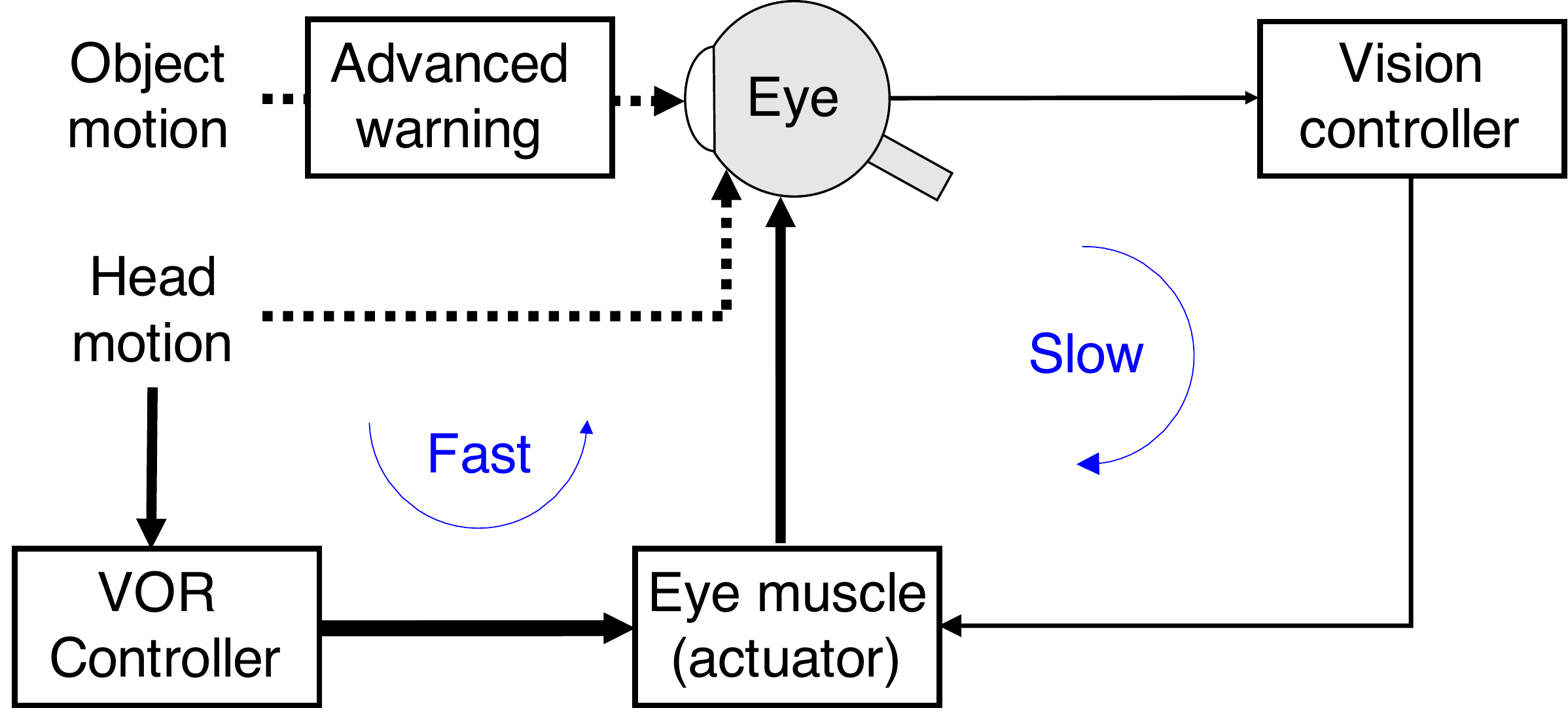}
         \label{fig:model visual-system}
\end{subfigure}
\begin{subfigure}[b]{\linewidth}
         \caption{Sensorimotor control system used in lateral control in trail following}
         \centering
	\includegraphics[width=0.8\linewidth]{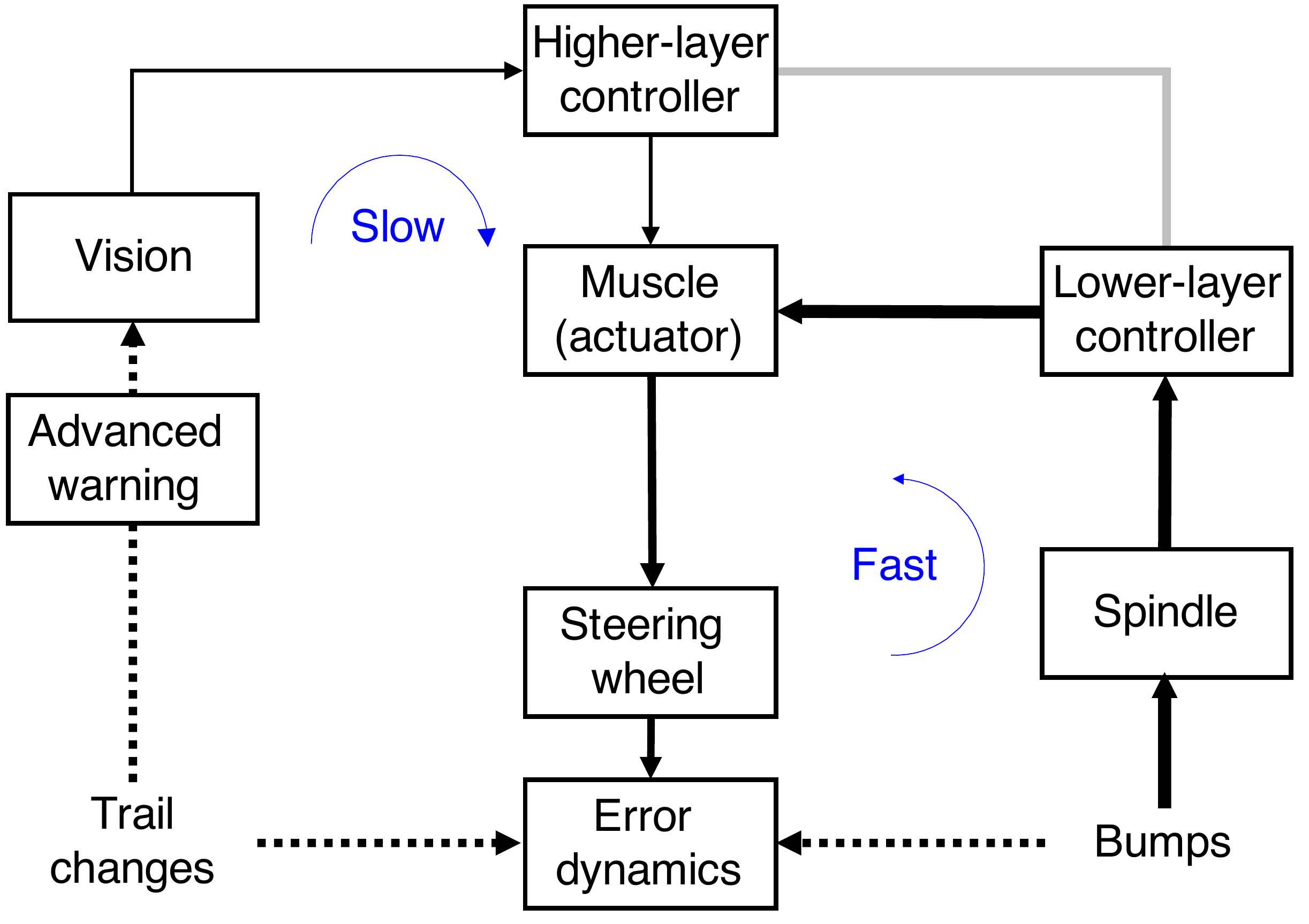}
	\label{fig:model biking-system}
\end{subfigure}

\caption{
\textbf{System diagrams of sensorimotor controls used for biking in a mountain trail: oculomotor control for visual tracking and lateral control in trail following.}
\textit{(A)}  Diagram of two major feedback loops involved in the eye movement: visual loop and vestibular-ocular reflex (VOR) loop. Objects are tracked using the slow visual loop, while head motion is compensated for by the much faster VOR loop. 
\textit{(B)} Diagram of the basic sensorimotor control model for our experiment that simulates lateral control in trail following.  
Each box is designated by its function: sensing and communication (\eg vision, muscle spindle sensor, vestibulo-ocular reflex), actuation (muscle), and computation (high-layer planning and tracking and low-layer reflexes and reactions). Depending on the hardware details, they may be quantized (discrete valued), have time delays, experience saturation, and be subject to noise. The trail ahead can be seen in advance, but the bumps and other disturbances are unanticipated. The line thickness indicates the relative speed of the pathway (thicker lines for faster pathways).
}
 \label{fig:model}
\end{figure}

In this study, we investigated the influence across component levels and integration across control layers in a psychophysical task related to mountain biking. 
Our study of sensorimotor integration is a first step toward bridging the persistent gap between the hardware limitations and systems performance. Our results suggest the importance of layering and diversity: The diversity between layers can be exploited to achieve {\em both} fast {\em and} accurate performance despite imperfect hardware that is slow or inaccurate. 

\subsection*{Control systems involved in mountain biking} Successfully riding a mountain bike down a bumpy, curved trail requires remarkable sensorimotor performance through the effective integration of many subsystems, including oculomotor control for planning, lateral control for trail following, and balancing in rough terrain. 

The oculomotor control system uses a layered control architecture to maintain fixation on a visual target while bouncing down a trail. The vestibulo-ocular reflex (VOR) compensates head jostling with fast feedforward circuits in the brainstem, and a slower feedback loop from visual cortex pursues moving targets (Fig.~\ref{fig:model visual-system})~\cite{lac1995learning,lisberger2010visual}. In addition, the cerebellum monitors proprioceptive inputs from muscles and efference copies from motor commands. This predictive feedback modulates the gain of the VOR in the context of the current state of the body and intended actions~\cite{coenen1995Vergence}.

The lateral position control for trail following uses an architecture with one layer that plans the trajectory and another layer that stabilizes against bumps and rocks on the ground (Fig.~\ref{fig:model biking-system}). Trajectory planning takes place in the cerebral cortex and basal ganglia using visual information of approaching obstacles such as trees and winding trails. The delay in visual processing between the retina and the eye muscles during smooth visual pursuit is around 100 milliseconds~\cite{Lisberger2015SmPur}. This higher layer of processing interacts with a lower layer having faster feedback loops in the spinal cord that control deviations from the desired trajectory generated by a bumpy road.

\begin{figure}[hb!] 
\begin{subfigure}[b]{\linewidth} 
\caption{Visual input to the subjects} \centering \includegraphics[width=\linewidth]{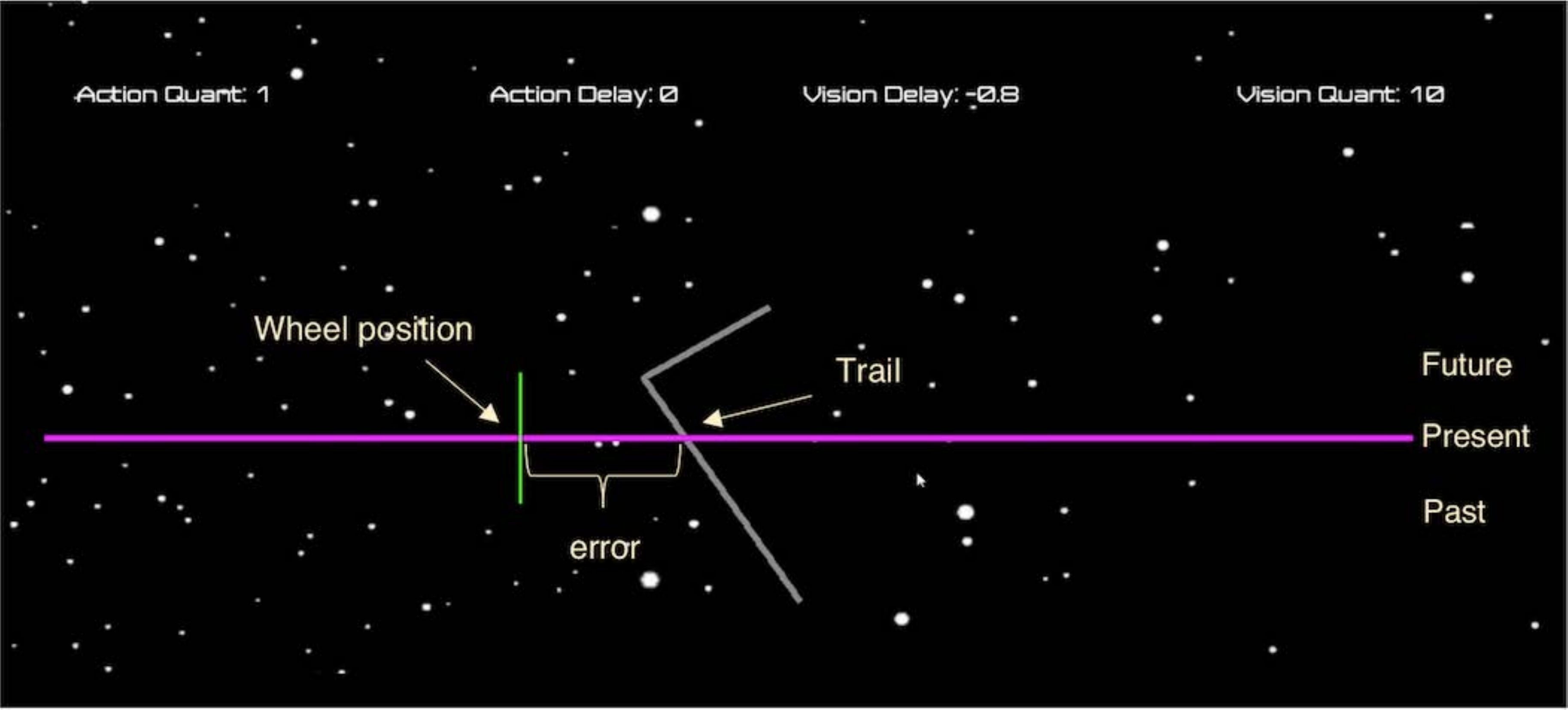} \label{fig:monitor} \end{subfigure} 
\begin{subfigure}[b]{\linewidth} 
\caption{The experimental setup} \centering \includegraphics[width=\linewidth]{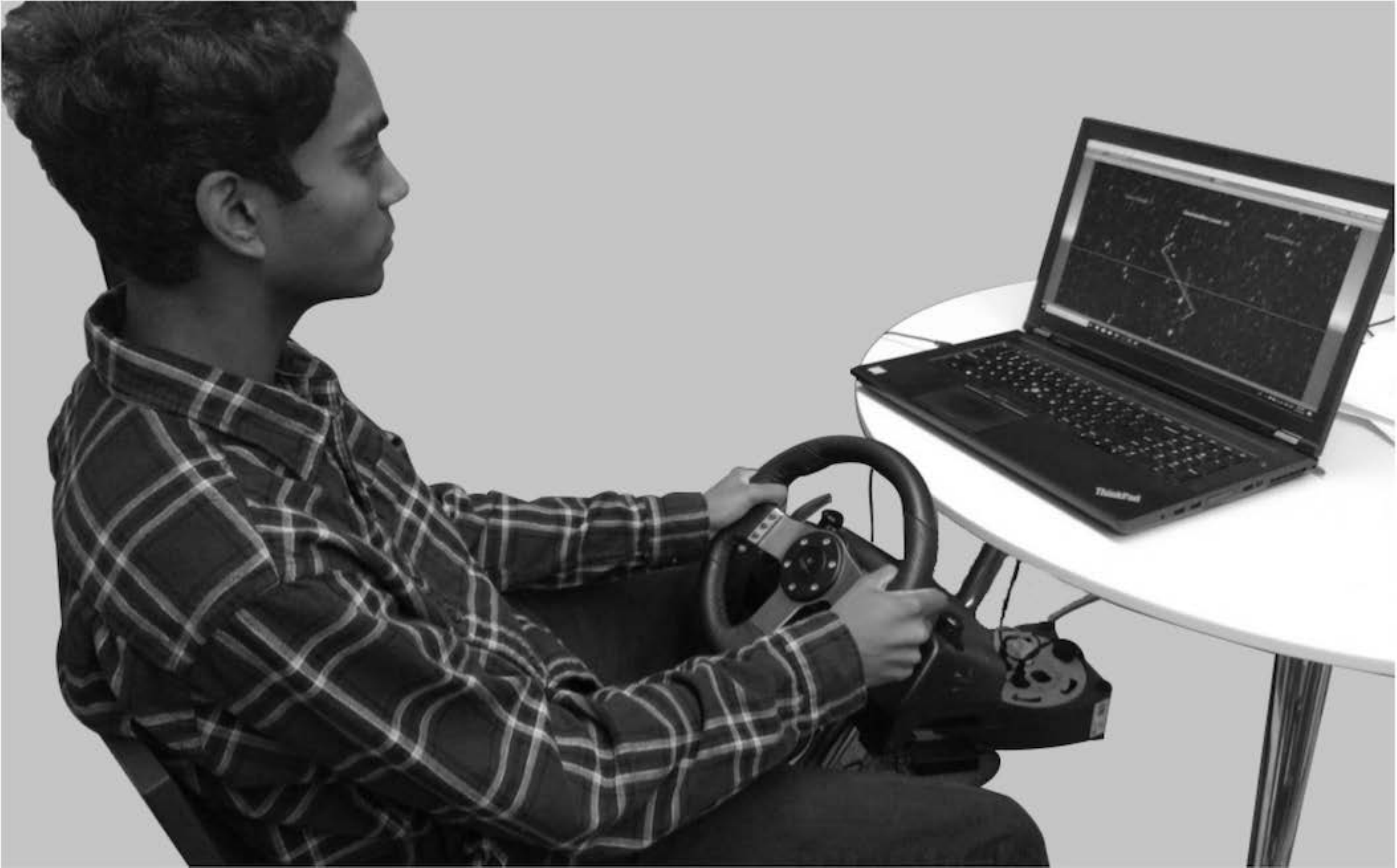} \label{fig:experiment} \end{subfigure} 
\caption{
\textbf{Video monitor interface for the biking task.}
\textit{(A)} Players see a winding trail scrolling down the screen at a fixed speed, and with a fixed advanced-warning (the visible trial ahead), both of which can be varied widely. The player aims to minimize the error between the desired trajectory and their actual position using a gaming steering wheel. \textit{(B)} Bumps are added using a motor torque in the wheel. Experiments can be done with bumps only or trails only, or both together, and with varying trail speed and/or advanced-warning, and with additional quantization and/or time delay in the map from wheel position to players' actual position.
}

\label{fig:interface}
\end{figure}

\begin{figure}[ht!]
\begin{subfigure}[h!]{\linewidth}

        \caption{Errors in the case of bump only, trail only, and both}
        \centering
\includegraphics[width=0.8\linewidth]{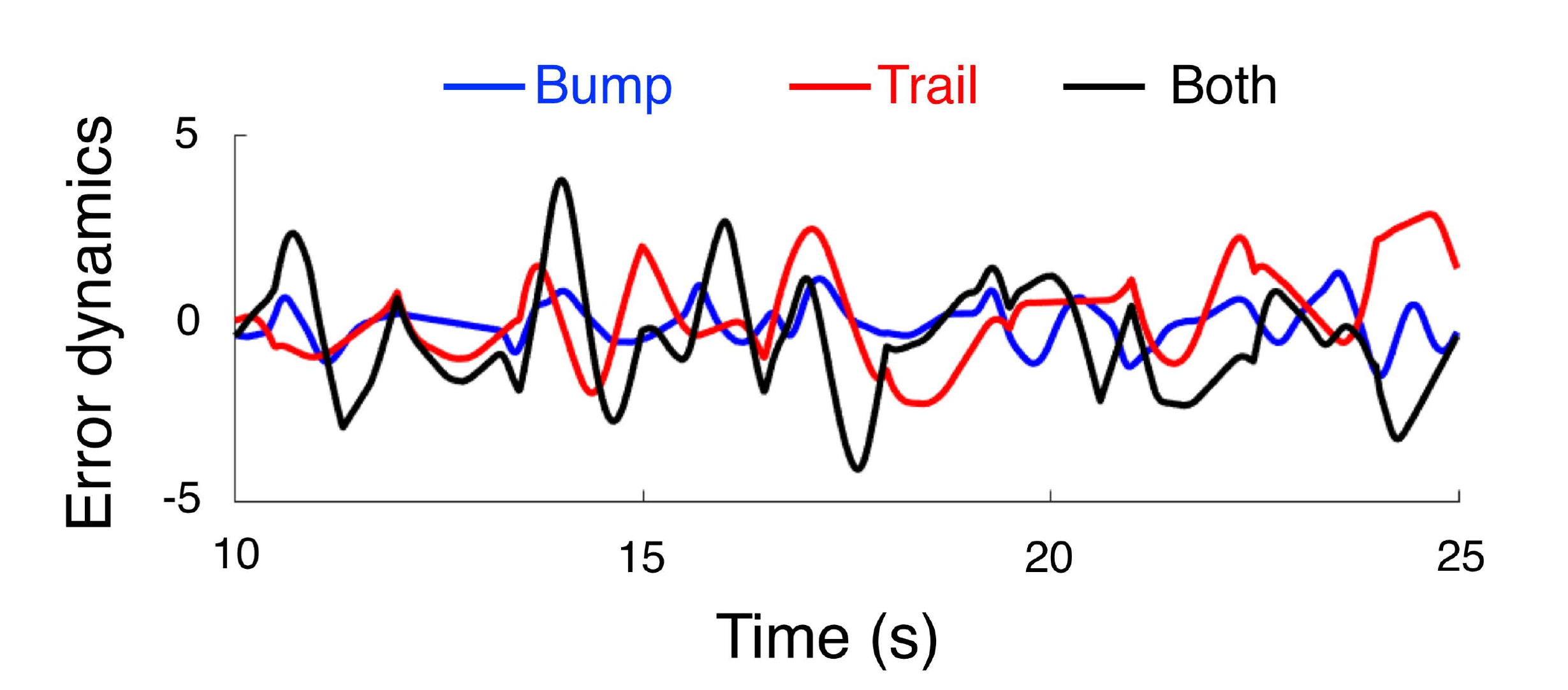}

\end{subfigure}

\begin{subfigure}[b]{\linewidth}
 \caption{Dynamics of the additive error in the cases of bump only and trail only versus error in the case of both bump and trail}
\center
\includegraphics[width=0.8\linewidth]{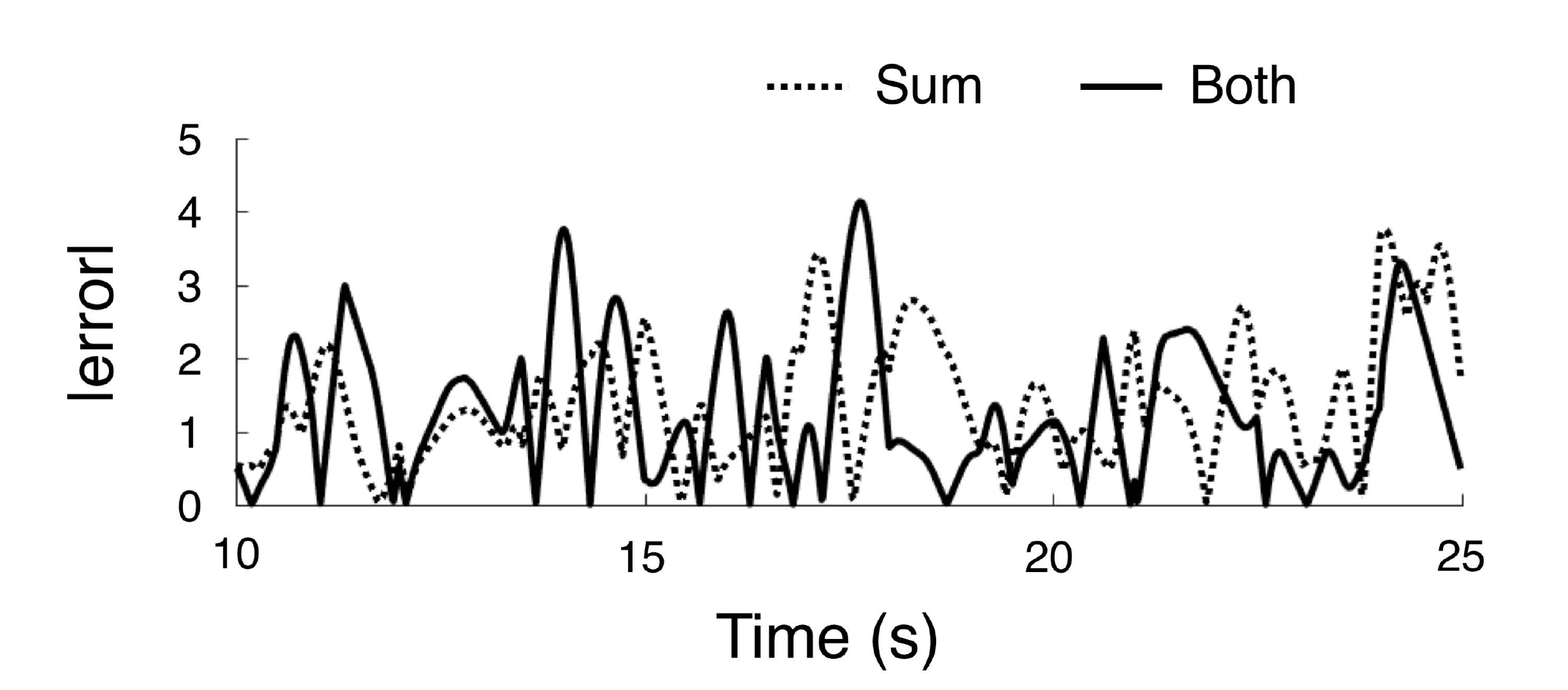}
\end{subfigure}

\begin{subfigure}[b]{\linewidth}
 \caption{Sizes of the additive error in the cases of bump only and trail only versus error in the case of both bump and trail}
\center
\includegraphics[width=0.8\linewidth]{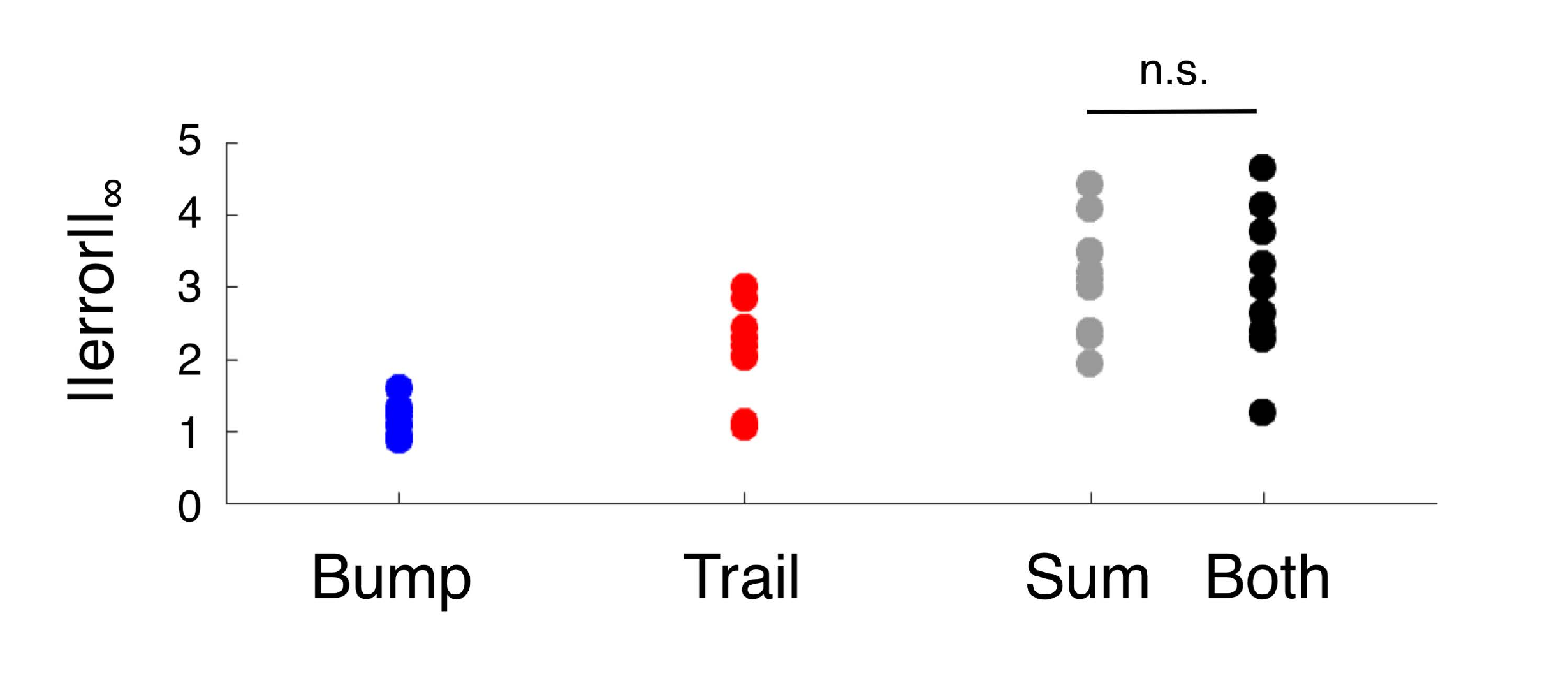}
\end{subfigure}

\caption{
\textbf{Total error and its decomposition into the error due to bumps and the error in tracking the trail.} 
\textit{(A)} Error dynamics from a task with only added bump, a task with only trail changes, and a task with both. \textit{(B)} The size of errors from the first two tasks and the error from the last task. \textit{(C)} Worst-case errors for the three cases and the sum of errors from the first two cases. Each dot denotes the worst-case error in 2 seconds. 
}
\label{fig:experiment-bumps-trail}
\end{figure}

\section*{Results}

\subsection*{Experiments}

\begin{figure}[ht!]

\begin{subfigure}[b]{\linewidth}

        \caption{Theoretical system SATs}
        \centering
	\includegraphics[width=0.9\linewidth]{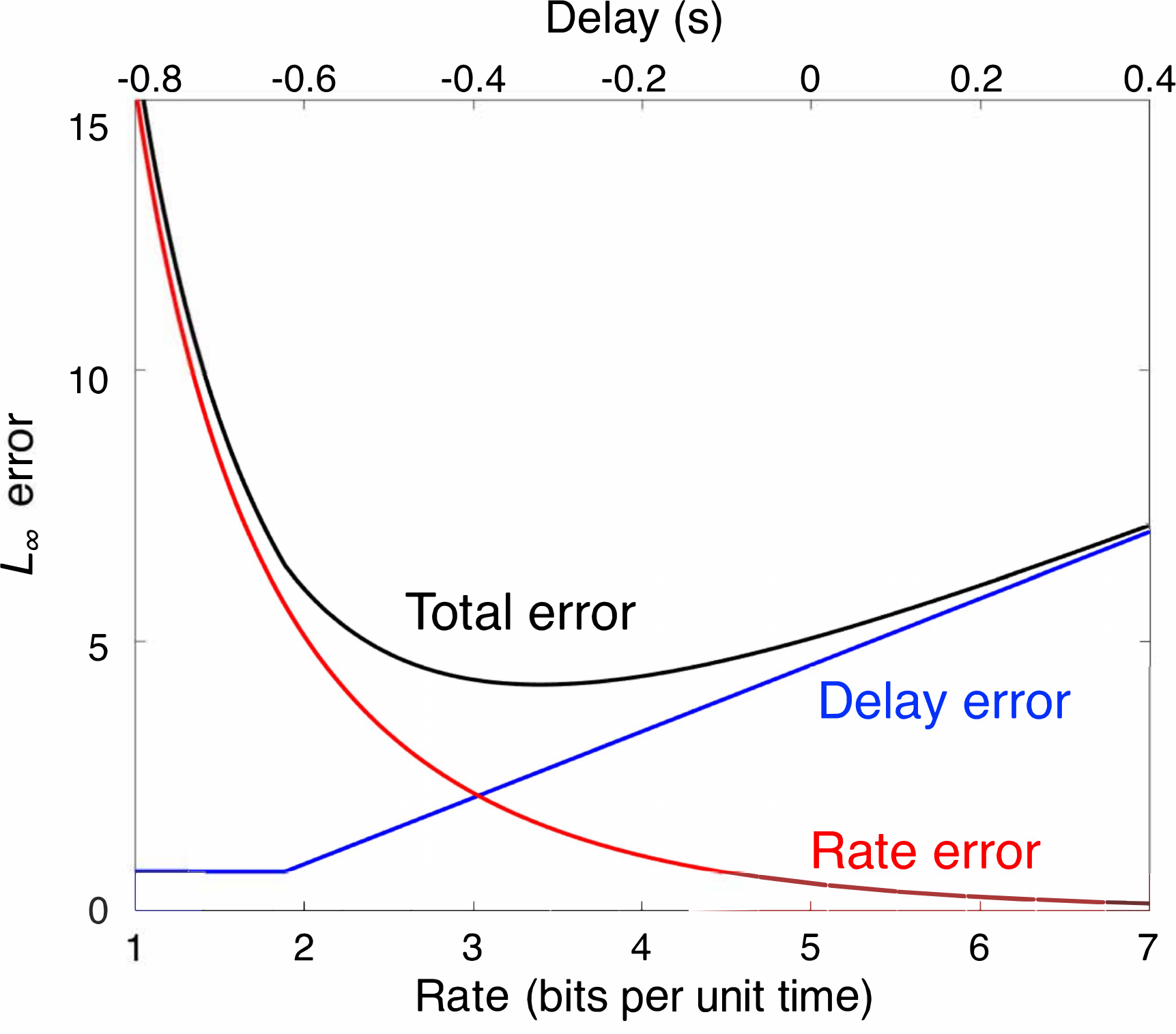}
	\label{fig:SAT-system-theory}
\end{subfigure}

\begin{subfigure}[b]{\linewidth}
	\caption{Experimental system SATs}
	\centering
	\includegraphics[width=0.9\linewidth]{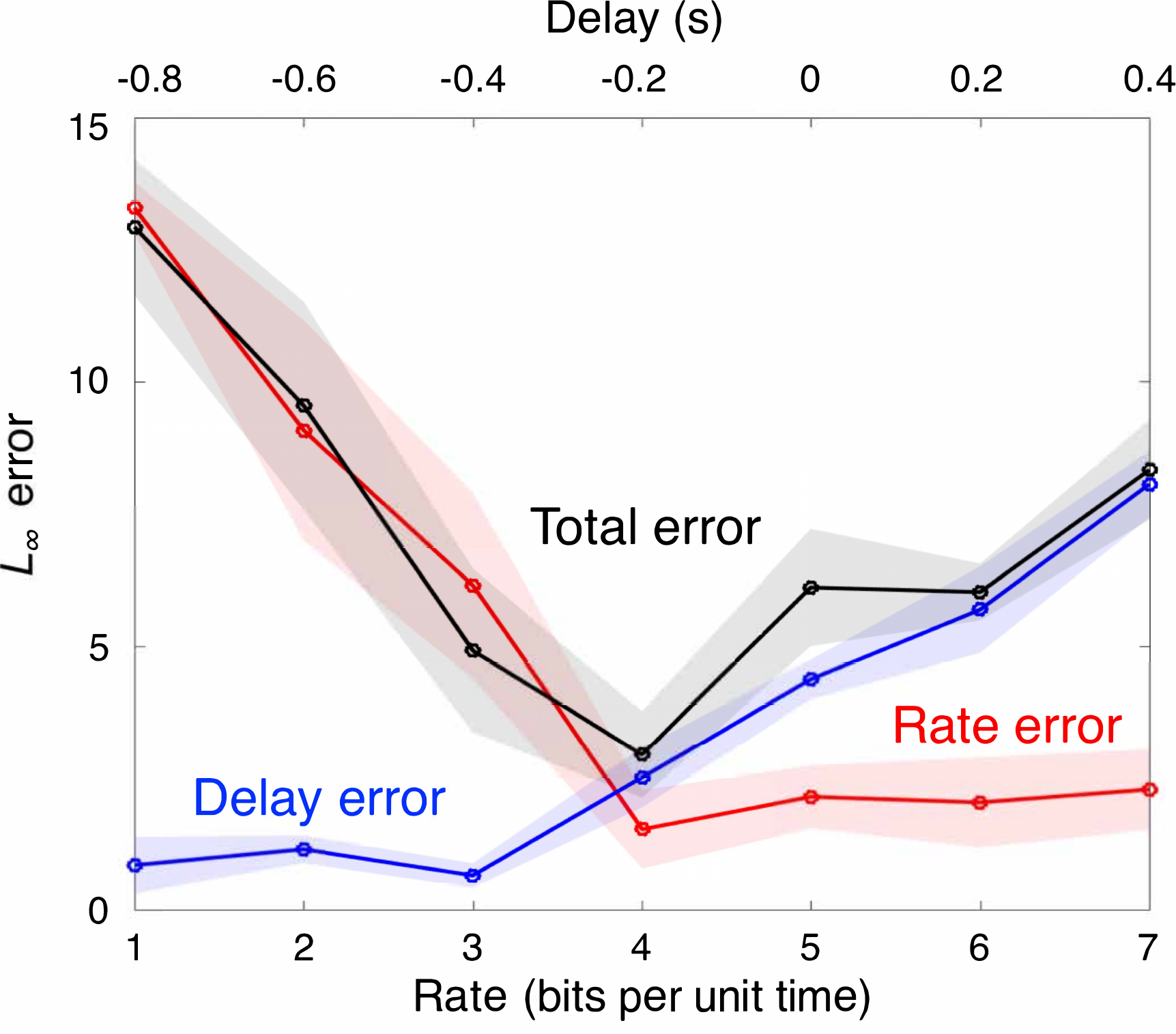}
	\label{fig:experiment-SAT}
\end{subfigure}
\caption{
\textbf{System SATs in the biking task.}
\textit{(A)} {\re Theoretical SATs. The delay error $\max(0,T)$ (blue), rate error $(2^{R}-1)^{-1}$ (red), and the total error $\max(0,T) + (2^{R}-1)^{-1}$ (black) in Eq.~\ref{eq:performance_det} are shown with varying component signaling delay $T_s$ and rate $R$ subject to the component SAT $T = (R-5)/20$.}
\textit{(B)} Empirical SATs. The error under an added delay (blue), the error under added quantization (red), and the error under added delayed plus quantization (black) are shown. In the last case, the added delay $T$ and quantization rate $R$ subject to the component SAT $T = (R-5)/20$. The dot shows the averaged error of 4 subjects, and the shadowed area indicates the standard error of the mean for these subjects. 
}
\label{fig:SAT}
\end{figure}

{\ed We developed experimental tasks and corresponding sensorimotor control models that mimicked three aspects of mountain biking: compensation by the spinal cord for the random shaking coming down the trail, the anticipation of turns in the trail by the visual system, and the stabilization of images on the retina by the oculomotor system to compensate bouncing. We experimented with the first two aspects using driving experiments and the last aspect using a few simple tasks that that the reader can easily perform. Many other aspects of biking are left out, but by focusing only on these aspects, we are able to make testable predictions. 

We performed two driving experiments: The first is to test the interactions between layers, and the second is to test the errors caused by delays and rate limits in control within a layer. In the two experiments, subjects follow the trail on a computer screen and control a cursor with a wheel to stay on the trail. The goal of the subjects is to minimize the errors between the desired and actual trajectories shown in a computer monitor by moving the steering wheel (Fig.~\ref{fig:interface}, see Materials and Methods for details). 

In the first experiment, the higher-layer and the lower-layer are coordinated. We compared how subjects' control behaviors and the resulting errors differ in three settings: 1) when there are random force disturbances to the steering wheel due to bumps on the ground (denote as 'Bump only'), 2) when the trail trajectory is curved and changes direction (denote as `Trail only'), and 3) when both exist (denote as `Both'). 
Rejection of bump disturbance in the first and last settings is likely to be performed at the lower layer reflex, while trajectory following in the second and last settings is likely to be performed at the higher layer planning. 

The experimental results are shown in
Fig.~\ref{fig:experiment-bumps-trail}. The observed error in setting 3 (with both bumps and trail curvature) positively correlated with the sum of the errors from the first two settings with either bumps or trail curvature (Pearson correlation coefficient $=0.57$), suggesting the two signals tended to have consistent sign and amplitude. Moreover, the two signals showed no significant difference in the two-side t-test analysis. 
The results suggest that the two layers could be analyzed separately. This separability motivates the modeling of each layer separately and to further decompose the errors into those caused by neural signaling delays or rate limits in the control loop. 

The impact of neurophysiological limits was studied in the second experiment. We observed changes in lateral control error in three settings: when external delays are added in the display, when external quantizers are added in the actuation effect of the steering wheel, and when both are added. These manipulations served as noninvasive probes for how component SATs constrain the system SATs. The lateral errors in the three settings are shown in Fig.~\ref{fig:experiment-SAT}, and their corresponding theoretical prediction is shown in Fig.~\ref{fig:SAT-system-theory} (see the modeling details in the next section). The bridge between the SATs at the two levels highlights the benefits of the heterogeneity observed in nerves  (Fig.~\ref{fig:axons}) and the advantages of layering in sensorimotor control (\eg Fig.~\ref{fig:model}). 

Our experiment primarily focused on the layers involved in lateral control. In both experiments, the head was relatively stable, and the errors of image stabilization on the retina by VOR, though essential in mountain biking, is negligible. Another important layer that was not included in the biking game was bike balance and turning, skills that must be learned before trail following. 
}

\begin{table}
\caption{Parameters in the basic model.}
\vspace{-5mm}
\begin{center}
\begin{tabular}{@{\extracolsep{5pt}}ll} 
\\[-1.8ex]\hline 
\hline \\[-1.8ex] 
\textbf{Parameter} & \textbf{Description} \\
\hline 
$x(t)$ & Error at time step $t$  \\
\hline
$\mathcal K$ & Controller  \\
\hline
$T_s \geq 0$ & Signaling delay  \\
\hline
$T_a \geq 0$ & Advanced warning  \\
\hline
$T_i \geq 0$ & Internal delay \\
\hline
$T = T_s + T_i - T_a$ & Total delay 
\\
\hline
$R$ & signaling rate (bits per unit time) 
\\
\hline
$\lambda$ & Cost associated with the resource use
\\
\hline
\end{tabular}

\centering
\label{tab:notations}
\end{center}
\end{table}

\subsection*{Connecting component and system SATs}

To connect the SATs between the two levels, we developed a robust control model that characterizes the system-level SATs imposed by component-level SATs and used the model to explain the experimental observations. We modeled the error dynamics between the  actual lateral position of the subjects and the center of the trail  as follows: 
\begin{equation}
\label{eq:plant}
x(t+1) = x(t) + w(t) +  u(t) , 
\end{equation}
where $x(t)$ is the lateral error, and Eq.~\ref{eq:plant} relates the future error $x(t+1)$ with the previous error $x(t)$, the uncertainty $w(t)$ (bumps and trail changes), and the control action $u(t)$. In the brain, the control action $u(t)$ is generated from many sources: 1) sensory information arising from visual inputs, proprioception from stretch receptors in muscles and acceleration from vestibular organs in the inner ear; 2) communication pathways through sensory and motor nerves; 3) computation in the central nervous system, including the spinal cord, cerebellum and cerebral cortex; and 4) actuation by muscles in the eyes and arms.

{\ed This simple model captures the bicycle dynamics and controller structures in the process and constraints that generate $u(t)$. This is not an all-encompassing model of all the biomechanics (\eg muscle mechanical properties, vesticulo-ocular reflex gain adjustment) and control loops (\eg physiological reflex loops) that are involved in mountain bike riding. Rather, it abstracts out the component delay and data rate, which are explained below, and consider the fundamental tradeoffs induced by these constraints in system performance. This abstraction allows us to focus on the mechanism to exploit diversity, which we believe is universal but has heretofore been ignored. Our approach, which is commonly used to tease apart the complexity of biological systems, does not deny the existence of the underlying complexity but will isolate each component from the complexity in order to nail down the scientific hypothesis worthy of further investigation. }


In the experiment, these are programmed by the software and can be made to be arbitrarily hard.

The feedback loop from sensor measurement $x(t)$ to control action $u(t)$ has a latency of $T_u := T_s + T_i$ with a signaling rate $R$, where $T_s$ models the nerve signaling delay, $T_i$ models other internal delays in the feedback control loop (including both sensory and motor delays), and $R$ is the maximum rate at which axons can transmit information. The feedforward loop from disturbance $w(t)$ to the control action $u(t)$ has an advanced warning of $T_a$. Advanced warning occurs when bikers view the future trail trajectory $T_a$ steps ahead, before it influences the error dynamics, which allows predictions to be made and muscle tone changes to occur ahead of time \cite{massion1992movement,bouisset2008posture}. The value of $T_a$ depends on the speed of the rider and the features on the trail. The delay from the moment the error dynamics are impacted by a disturbance to the moment the control acts against the disturbance is the latency minus warning, \ie $T := T_u - T_a = T_s + T_i - T_a$. The list of notations is shown in Table~\ref{tab:notations}.

\textit{Component SATs.} Next, we characterize the trade-off between nerve signaling delay and rate limit arising from the fixed spatial and metabolic cost to build and maintain axons~\cite{perge2009optic,perge2012axons,sterling2015principles,laughlin2003communication}. Specifically, nerves with the same cross-sectional area can either contain many small axons or a few large axons (Fig.~\ref{fig:axons}), which inevitably leads to SATs in neural signaling. The specific forms of SATs depend on how the nerves encode information~\cite{mainen1995reliability,salinas2001correlated,bodelon2005temporal,fox2010encoding}, and a wide range of time-based and rate-based codes are found throughout brains~\cite{Rieke1997Spikes}. 

In the spike-based encoding scheme, information is encoded in the presence or absence of a spike within each time interval, analogous to digital packet-switching networks~\cite{salinas2001correlated,srivastava2017motor}. For example, spike-based coding is found in many subcortical structures, such as spatial localization in the auditory system encoded as time delay between the two ears~\cite{Grothe2010Auditory}, and spike timing in the cerebral cortex regulates synaptic plasticity~\cite{tiesinga2008regulation}. It has also been observed that many types of neurons can generate spikes with accurate timing~\cite{mainen1995reliability,fox2010encoding}, which is typically required in spike-based or time-based codes.

Assuming all axons have the same size, the component SATs can be shown to satisfy 
\begin{align}
\label{eq:delay-rate-spike1}
&R = \lambda T_s  ,
\end{align}
{\ed 
where $R>0$ is measured by bits per unit time, $T_s>0$ is measured in unit time, and $R$ and $T_s$ should use an identical time unit.} The constant $\lambda(>0)$ is proportional to the spatial and metabolic resources required to build and maitain the axons. In rate-based encoding, the SATs is approximated using the information capacity of a communication channel of Poisson type: 
\begin{align}
\label{eq:spike-rate-code}
 &R =  \frac{1}{2} \lambda  T,
\end{align}
where $R$, $T$, and $\lambda$ are
the same variables as in Eq.~\ref{eq:delay-rate-spike1}. The value of $\lambda$ depends on the nerve. For example, proprioceptive nerves often have large $\lambda$, which allows lower latency with higher data rate compared with unmyelinated pain fibers, which are slow but have a high signaling rate. For a fixed resource level $\lambda$, the same rate $R$ can be achieved with half of the delay using spike-based encoding than with (less efficient) rate-based encoding. 

{\ed Eqs.~\ref{eq:delay-rate-spike1} and~\ref{eq:spike-rate-code} characterize the amount of information that can be transmitted within the latency requirements for control and are derived as follows. First, the space and energy to build and maintain nerves, quantified by $\lambda$, are translated into the size vs. number tradeoff for axons. Next, the size vs. number tradeoff is converted into the latency and rate tradeoffs. Here, the speed at which an action potential travels depends on the axon size and the maximum firing rate depends on the metabolic energy that is available. These constraints, together with assumptions on how the information is encoded, determine the maximum signaling rates. A detailed derivation is given in the supplementary material. Our approach is different from the approach that uses asymptotic information theory to characterize the amount of information that can be transmitted without considering latencies.}

Note that the system performance limits of our model do \textit{not} require the component SATs to take the forms given in Eqs.~\ref{eq:delay-rate-spike1} and~\ref{eq:spike-rate-code}. Component SATs differ by encoding schemes and the presence or absence of myelination, noise and redundancy, and cross-talk between axons. Although we use the SATs in Eq.~\ref{eq:delay-rate-spike1} here, similar analysis can be performed for other component SATs.


\textit{System SATs.} When performing sensorimotor control, the component-level SATs constrain the system-level SATs. To characterize this relation, we first use robust control tools to find the errors as functions of the component-level signaling delays and rates in both deterministic and stochastic settings. 

The worst-case framework is suitable for modeling risk-averse sensorimotor behaviors, such as riding a mountain bike on a trail in the presence of the life-threatening  uncertainty is l~\cite{whittle1990risk,nagai1996bellman,sanger2010neuro,sanger2014risk}. {\ed When the disturbance is bounded in infinity norm, the worst-case error normalized by the size of the disturbance satisfies
\begin{equation} 
\begin{aligned}
\label{eq:performance_det}
\sup_{ \|  w \|_\infty \leq 1} \|  x\|_\infty \geq \max(0, T) + \left(2^{R} - 1 \right)^{-1}.
  \end{aligned}
\end{equation} 
This error $\sup_{ \|  w \|_\infty \leq 1} \|  x\|_\infty  = \sup \{  \|x\|_\infty / \|w\|_\infty\}$ captures the ratio of amplification or attenuation in worst-case error per unit size disturbance in worst case. This ratio can be used with different units. For example, if the sampling interval for control is $\tau$ seconds, then $T$ has a unit of sampling intervals, and $R$ has a unit of bits per sampling interval. If the disturbance has a size of $W$ cm per second, then the disturbance has size $|w(t)|\leq \tau W$ cm at each sampling interval, and the error is bounded by $\max(0, T) + \left(2^{R} - 1 \right)^{-1}$. Eq.~\ref{eq:performance_det} also applies when there is feedback, when the controller senses $x$, and  feedforward control, when the controller senses $w$. 
}

{\ed The average-case framework is more applicable to risk-neutral sensorimotor behaviors, such as riding a mountain bike across a broad field, where fatal risk is minimal~\cite{franklin2011computational,todorov2002optimal}. The precise formulations of the control problem for both cases are given in the supplementary material. When the disturbance has zero mean and bounded variance, the steady-state mean squared error normalized by error variance satisfies 
\begin{equation}
\label{eq:performance_stc}
\begin{aligned}
\sup_{ \mathbb E[w] = 0, var(w) = 1} \mathbb E[  x^2 ]  \geq \max(0, T) + \left(2^{2R} - 1 \right)^{-1}.
\end{aligned}
\end{equation}
This error $\mathbb E[  x^2 ] =  \mathbb E[  w^2 ] $

This error captures the ratio of amplification or attenuation in average error per unit variance in disturbance. } 

{\re The error bounds in both cases (Eq.~\ref{eq:performance_det}--\ref{eq:performance_stc}) are qualitatively similar: both bounds decompose into two terms. The shared first term, $ \max(0, T)$, only depends on the total delay and thus can be considered as the delay error. The other terms, $(2^{R} - 1)^{-1}$ and $(2^{2R} - 1)^{-1}$, depend only on the signaling rate and can be considered the rate error. Here, the units of the delay and rate errors are based on control, which are measures of system performance, rather than time or information measures (\eg bits), which are the units used in the signaling delay and rate at the component level.

This decomposition of errors is consistent with the experimental observation that the error for the trials with both added delay and added quantization was approximately the sum of the errors for the trials with the delay and the quantization added separately (Fig~\ref{fig:experiment-SAT}, see Materials and Methods for details).} {\ed The delay error, the rate error, and the total error in the experiment contain the internal errors of the subjects' sensorimotor control system in addition to the error caused by added delay, the error caused by added quantization, and the error caused by added delay and quantization, respectively. Therefore, the total error equals the sum of delay error (error due to added delay plus internal error) and rate error (error due to added quantization plus internal error) minus internal error. When the error due to added delay and quantization vanishes to zero, the delay and rate error approximately equals the internal error, and the total error converges to the delay and rate error because the delay and rate error (approximately the internal error) plus the rate and delay error minus the internal error (approximately the delay and rate error).}  
Beyond the worst case framework described above, the same conclusion holds for the stochastic setting (experiment and results are in section 4 of the supplementary material).

We are now ready to characterize how the SATs at the component level impact the SATs at the system level. Combining the component SATs in Eq.~\ref{eq:delay-rate-spike1} for spike-based encoding, the error bound of Eq.~\ref{eq:performance_det} in the worst case, and $T = T_s + T_i - T_a$, we obtain the system SATs (the influence of the neural signaling constraints on sensorimotor control) in Fig.~\ref{fig:SAT-system-theory}. Here, a similar analysis can also be performed with other forms of component SATs (encoding schemes) or with the error bound of Eq.~\ref{eq:performance_stc} in the average case. 

Increasing the delay in the feedback loop increases the delay errors, while increasing the rate leads to a large decrease in the rate errors. Thus, delays can cause small disturbances to escalate into larger errors, and increasing the rate reduces errors exponentially in the context of control. These properties of the fundamental limits hold for rate constraints imposed on the sensing, communication, and actuation units. Intuitively, actuation quantization gives errors in control, whereas sensor and communication quantization gives errors in state estimate, which in turn leads to equal-sized errors in control. However, this condition may not hold if feedforward/predictive control is used to compensate for the same set of disturbances.

The minimum error is achieved when the deleterious effects of the nerve signaling delay and inaccuracy are both controlled within a moderate range. Thus, both the nerve composition that minimizes the delay of nerve signaling and the composition that maximizes its rate work together, resulting in suboptimal performance. {\re In particular, choosing components that optimize the signaling rate, which is often done in models based on asymptotic information theory, may lead to large delays and less robust sensorimotor control. When resources are limited, optimization  must balance the impacts of both signaling delay and signaling rate. The consequences of this trade-off are explored in the Discussion.}


\begin{figure}[ht!]
\center

\begin{subfigure}[b]{\linewidth}
	\caption{Optimal worst-case errors and their compositions for varying net delay and net warning}
	\centering
\includegraphics[width=\linewidth]{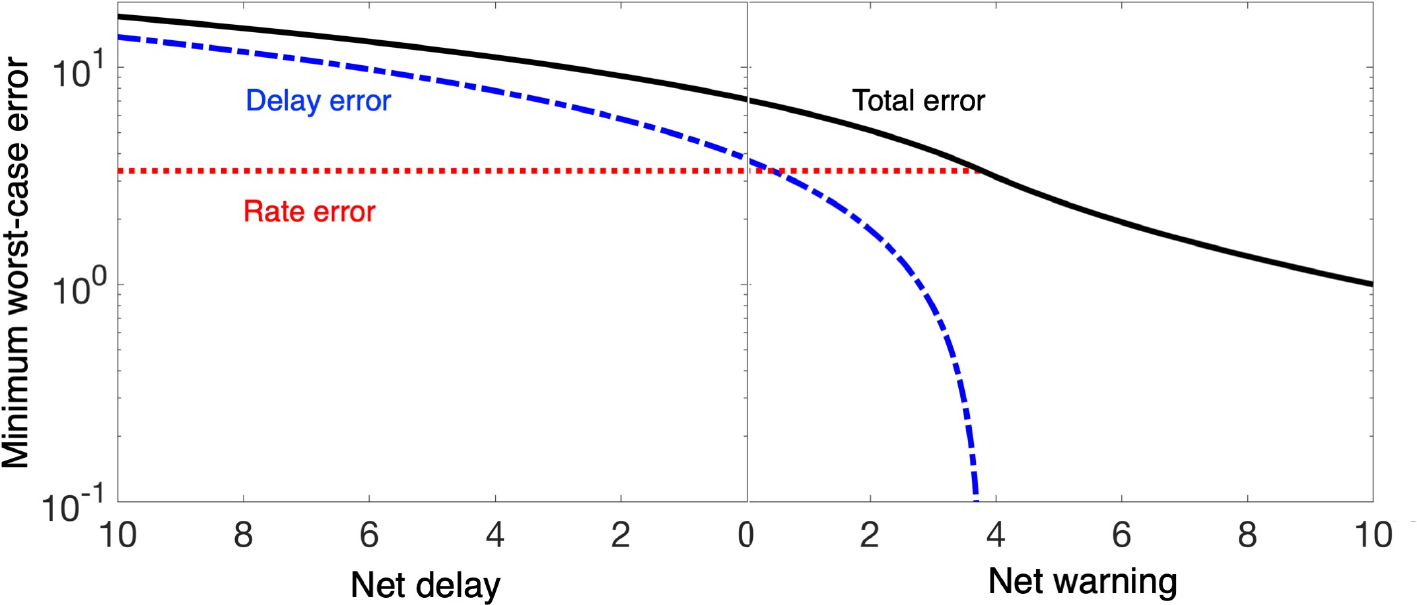}

\end{subfigure}

\begin{subfigure}[b]{\linewidth}
	\caption{Optimal total delays, signaling delays and rates for varying net delay and net warning}
	\centering
\includegraphics[width=\linewidth]{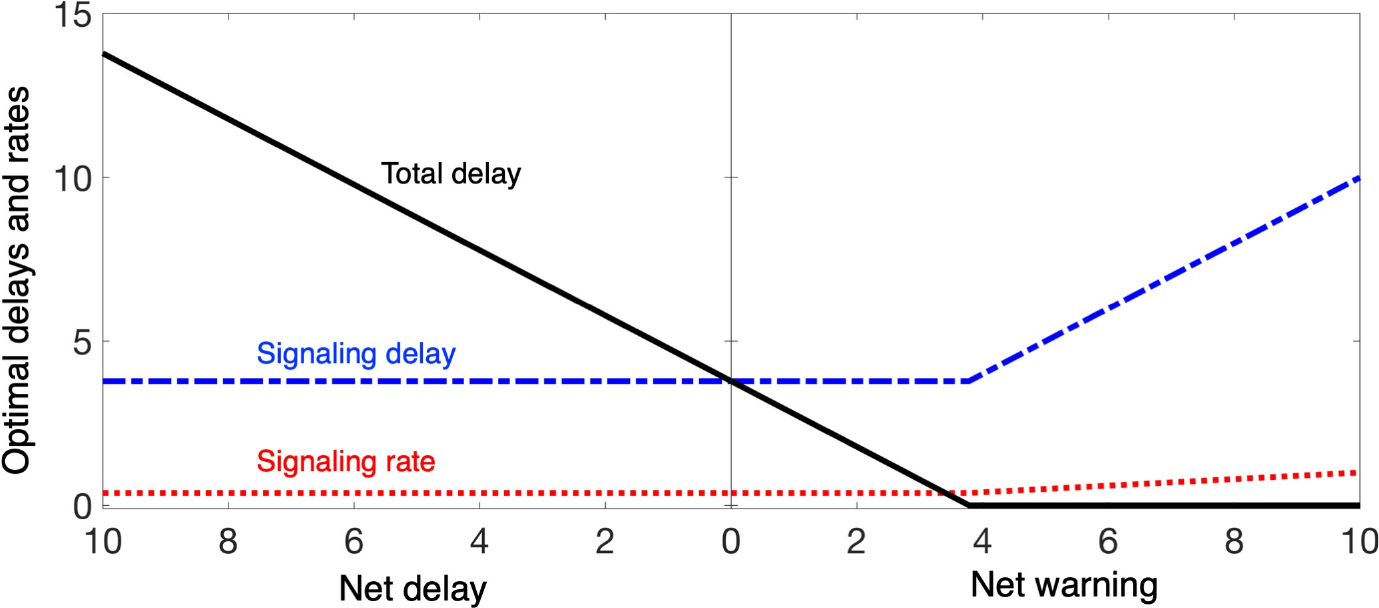}

\end{subfigure}

\caption{\textbf{Delayed reaction versus advanced planning for visual tracking.} \textit{(A)} The minimum total error Eq.~\ref{eq:performance_det} subject to the component SATs in Eq.~\ref{eq:delay-rate-spike1} and its composition (the delay error and the rate error) are shown for varying net delay $T_i - T_a (\geq 0)$ or net warning $T_a - T_i (\geq 0)$. \textit{(B)} The optimal signaling delay $T_s$, total delay $T (:= T_s + T_i - T_a)$, and rate $R$ for varying net delay or net warning. In both A and B, the resource to build and maintain axons are assumed to be fixed and are set to  $\lambda = 0.1$. 
}
\label{fig:SAT-delay-vs-advanced}
\end{figure}


\subsection*{Oculomotor control for visual tracking}
We now apply our theory to the layered control architectures for visual object tracking. {\ed The visual tracking of a moving object by smooth pursuit involves two major control loops: a fast feedforward VOR loop that compensates for head motion and a slower visual feedback loop through the visual cortex  (Fig.~\ref{fig:model biking-system})~\cite{lac1995learning,lisberger2010visual}.}

The vestibular inputs project to both the vestibular nucleus and the cerebellar cortex, which in turn projects back to the vestibular nucleus. This feedback loop from the cerebellar cortex is important for tuning the gain of the feedforward pathway to the vestibular nucleus. This tuning allows adaptation to the growth in head size during development and optical gain changes from new eyeglasses. The cerebellar inputs also correct the gain for changes in fixation distance and torsional head movements~\cite{Snyder1992, coenen1995Vergence}. Drifts across the retina due to unmatched gains are compensated by the visual system, which maintains vernier hyperacuity in the arcsecond range (2\% of the diameter of a cone photoreceptor in the fovea) for images drifting up to 3$^\circ$/sec~\cite{Westheimer1975hyperacuity}.

From a control perspective, an important difference between the two loops is their levels of advanced warning. The VOR loop reacts rapidly after the head moves. We call this regime the \textit{delayed reaction} of the VOR loop, in which the uncertainty $w(t)$ becomes accessible to the controller \textit{after} $w(t)$ affects the error dynamics, giving rise to positive net delay $T_i - T_a \geq 0$. In biking, vision allows looking ahead down the trail, which translates into a net advanced warning with enough look ahead. But in VOR, this doesn't happen.  

In contrast, changes in the visual environment are highly  predictable, so the visual loop can plan eye movements in advance, a  negative net delay.   
We call this regime the \textit{advanced planning} of the visual loop, in which the uncertainty $w(t)$ becomes accessible to the controller \textit{before} $w(t)$ affects the error dynamics, giving rise to positive net warning $T_a - T_i \geq 0$ (negative net delay). These two regimes are qualitatively different in their optimal choice of $T_s$ and $R$ for achieving the optimal robust performance, as illustrated in Fig.~\ref{fig:SAT-delay-vs-advanced}. 

\textit{(i) Delayed reaction:} When the net delay $T_i - T_a>0$ is large, the total error can be much larger than the size of the uncertainty $\| w \|_\infty$ and goes to infinity as $T_i  \rightarrow \infty$. This large error amplification is consistent with the all-too-familiar observation that even a small bump on a trail can cause a cyclist to lose control of the bike and crash. As $T_i$ increases, the delay error increasingly dominates the total error. Since the delay error largely contributes to the total error, the total error is minimized when $T_s$ is set to be small in return for small $R$. Therefore, a feedback loop in this regime performs better when it is built from a few large axons. Interestingly, the flat optimal delay/rate within the delayed reaction regime suggests that optimal performance can be achieved using one type of nerve composition for a broad range of advanced warnings. This property is beneficial because the net delay (defined from advanced warning) differs across different sensorimotor tasks. 

\textit{(ii) Advanced planning:} When the net warning $T_a -T_i > 0$ is large, the total error approaches zero as $R\rightarrow \infty$. This large disturbance attenuation is consistent with the observation that a cyclist can avoid obstacles given enough time to plan a response, such as taking a path around them or bracing against their impact. Given sufficiently large advanced warning $T_a$, the rate error increasingly dominates the total error because the growth in $T_s$ incurs no additional delay error. Since the rate error contributes largely to the total error, the total error is minimized when the signaling rate $R$ is set to be large at the expense of large signaling delay $T_s$. Therefore, a feedback loop in this regime performs better when it is built from many small axons. 


This prediction is qualitatively consistent with the anatomy of the human oculomotor system (Fig.~\ref{fig:axons}). The vestibular nerve, which transmits three-dimensional velocity information from the inner ear to the vestibular nucleus in the brainstem,  has $20,000$ axons with a mean diameter $3 \, \mu m$ and coefficient of variation $0.4 \, \mu m$.  These fast axons allow feedforward eye muscle control with a delay of approximately $10$ ms delay~\cite{bodelon2005temporal}. In contrast, the optic nerve carrying visual signals from the retina has approximately $1$ million axons with a mean diameter  $0.6 \, \mu m$ and coefficient of variation $0.5 \, \mu m$, significantly smaller but more numerous and with greater variability \cite{sterling2015principles,laughlin2003communication}. The optic nerve projects to the cortex through the thalamus, where visual signals are sequentially processed in several cortical areas before projecting back to subcortical structures that control eye movements.  As a consequence of this long loop, the  visual feedback delay is approximately $100$  ms. The reader is invited to perform an experiment to illustrate the consequences of layering in the oculomotor system.

\begin{figure}[ht!]
\center

\begin{subfigure}[b]{\linewidth}
\center
	\caption{Optimal delays and rates for varying advanced warning}
	\centering
\includegraphics[width=0.7\linewidth]{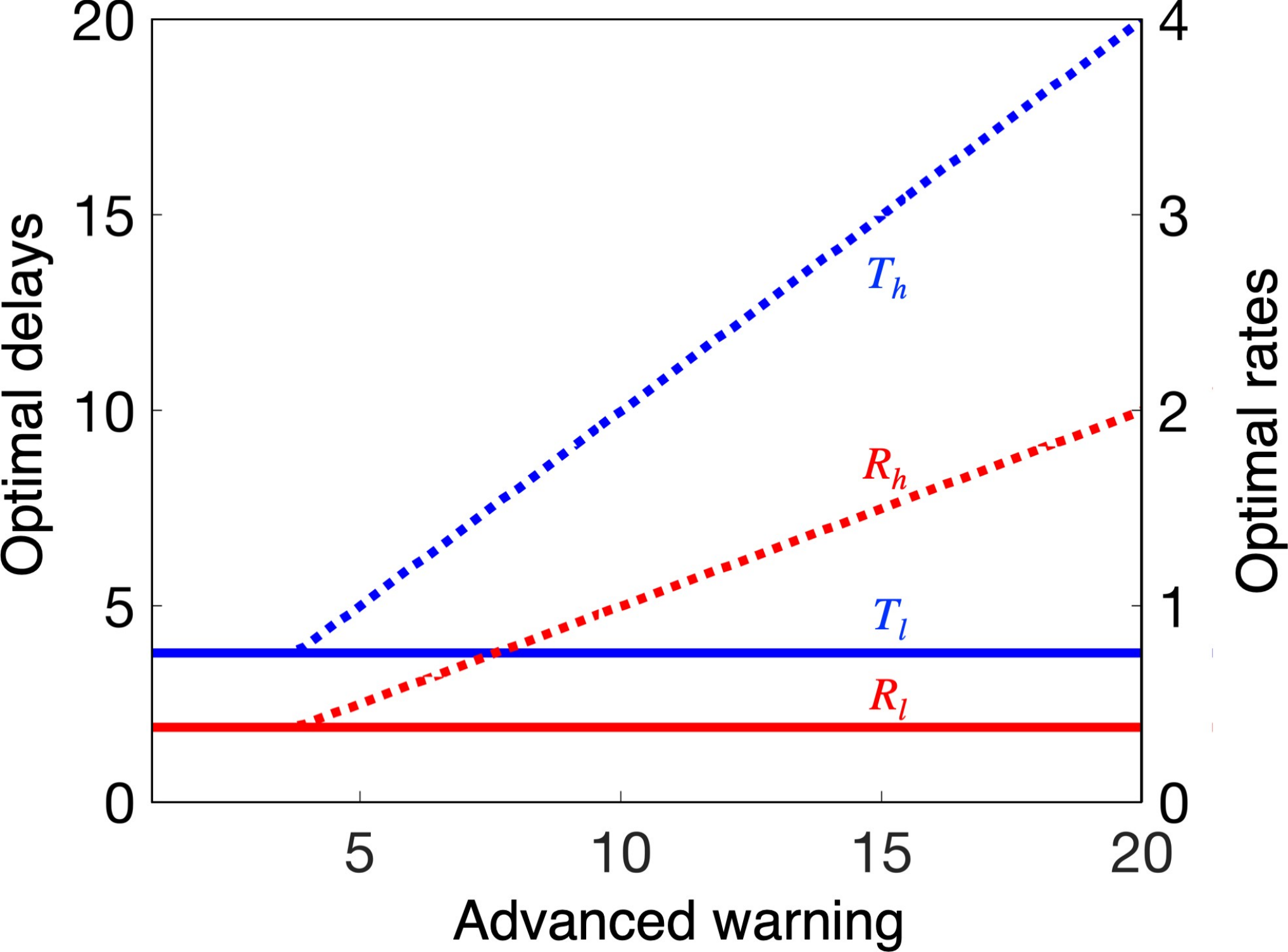}
\label{fig:loop-hete-joint-a}
\end{subfigure}

\begin{subfigure}[b]{\linewidth}
\center
	\caption{Optimal worst-case errors for varying advanced warning}
\includegraphics[width=0.7\linewidth]{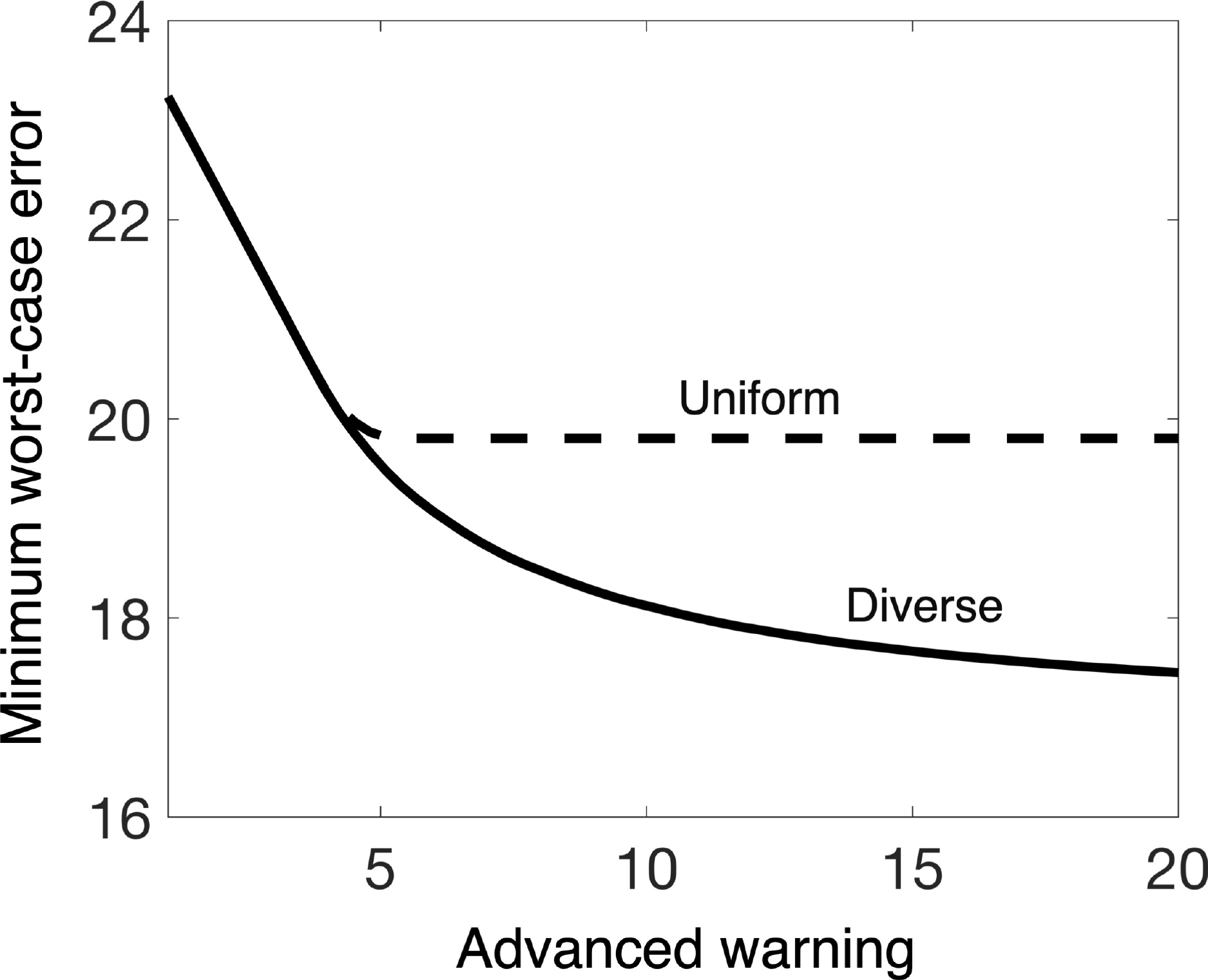}
\label{fig:loop-hete-joint-b}
\end{subfigure}

\label{fig:loop-hete-joint}
\caption{
\textbf{Planning and reflex layers for lateral control in trail following.} \textit{(A)} The optimal delays and rates for system performance Eq.~\ref{eq:optimal1} subject to the component SATs of the reflex layer $R_h = \lambda_h T_h$ where the delay in the reflex layer $T_i=10$, and the planning layer $R_\ell = \lambda_{\ell} T_{\ell}$ when the levels of advance warning $T_a$ are varied from 1 to 20. We set $\lambda_h=\lambda_{\ell}=0.1$ for the two layers.
\textit{(B)} Minimum error Eq.~\ref{eq:optimal1} for the case when the high-layer and low-layer are allowed to use components with diverse delay and rate or uniform delay and rate. The delays and rates in the diverse case do not have to be identical for both layers, whereas they are constrained to be identical in the uniform case, \ie $R_\ell = R_h$ and  $T_\ell = T_h$. In the diverse case, the high-layer controller can better exploit the advanced warning to minimize errors than in the uniform case. }

\label{fig:SAT-system}
\end{figure}

Layering diversity can be observed by tracking your hand moving left and right across the visual field with increasing frequency while holding the head still and comparing this with shaking the head back and forth (in a 'no' pattern) at an increasing frequency while holding the hand still. The hand starts to blur due to delays in visual object tracking at around 1-2 Hertz, whereas blurring due to the inability to compensate for fast head motion occurs at frequencies above 20 Hz. The difference is that the visual loop has lower levels of tolerable delays than the VOR loop. However, though slower, the visual loop is more accurate. 

Although both the VOR and visual layers have hard limits in speeds and accuracies individually, the limits do not translate into ‘inaccurate’ or 'slow' control at the system level because each layer is design to exploit the structures of the process to be implemented in that layer. Specifically, the limited signaling rate of the VOR loop does not compromise stabilization against head motions because this loop only requires three-dimensional velocity information in control. The visual feedback loop exploits the predictability in the visual environment to mitigate the latency in visual information processing. A separation of these two loops allows stabilization of head motion with a lower-dimensional velocity signal and visual object motion typically provides a large advanced warning. 

{\re The benefits of diversity between layers are visualized in Fig.~\ref{fig:DESS-visual}, which compares the system SATs when the VOR and visual layers use diverse delays and rates with a case when the delays and rates are uniform. Given the same amount of resources to build and maintain axons, the performance is more robust in the diverse case. }

This case study suggests an inaccurate but fast visual tracking layer and a slower but accurate VOR layer jointly create a virtual eye controller which is both fast and accurate. Although the component SAT imposes system trade-offs between minimizing the delay errors or rate errors in sensorimotor control, diversity deconstrains severe system SATs by using a slower but more accurate higher visual layer to reduce the rate cost and an inaccurate but faster lower reflex layer to reduce the delay cost. We call this \textit{diversity-enabled sweet spots} (DESSs): \ie the diversity between different layers helps achieve \textit{both} fast \textit{and} accurate sensorimotor control despite the slowness or inaccuracy of individual layers. {\re Sweet spots are in regions near the origin where delay errors and rate errors are both minimized.}

{\ed There are others layers in the oculomotor system.  For rapid saccadic responses to planned targets, the location from retinal sensors project directly to the superior colliculus, from which saccades are launched. For targeting a moving stimulus, peripheral retinal inputs have lower latencies than in the fovea.}

\subsection*{Visual and vestibular feedbacks for balancing} 
{\ed
Analogous DESSs can also be observed in the control processes used to balance unstable biking dynamics. Balancing uses a layered architecture involving visual, vestibular, and proprioceptive control loops. The development of the control system for balance begins in children 6 to 18 months old and is further enhanced with practice for more complex tasks such as biking. 
Visual, vestibular, and proprioceptive loops have diverse speeds and information rates, which complement each other to produce robust performance in balancing. Balancing with one leg is easy with normal visual and vestibular systems, and significant loss of balance with eyes closed often indicates proprioceptive or cerebellar injury. Standing on one leg is also harder with eyes closed than with eyes open because the vestibular loop without vision does not have access to the accurate information from the visual loop. Spinning or alcohol (or drugs) temporarily disrupt the vestibular control and increase the difficulty of standing in one leg. Unilateral or bilateral vestibular loss is also known to compromise the robustness of balancing and posture control~\cite{creath2002limited,black1989effects}.   
 }

\subsection*{Lateral control in trail following}

{\ed DESSs can also be observed in the layered control architecture used for lateral control in mountain biking. Planning loops at a high layer of visual processing in the cortex and basal ganglia track the trail. Spinal feedforward control compensates for large bumps and feedback compensates for small bumps, disturbances that are difficult to see. Below and above these two layers, a lower layer regulates muscle stiffness in anticipation of future bumps, and higher layers make cognitive decisions that are strategic. Here, we focus on the visual planning and reflex layers in the context of robust control and component diversity.}

To understand this mechanism, we use biking experiments (Fig.~\ref{fig:interface}) to simulate the lateral control in the mountain biking when the impact due to head and body movements are negligible. The lateral error dynamics is given by Eq.~\ref{eq:plant}, where $x(t)$ is the error, $w(t)$ is the disturbance, and $u(t)$ is the control action. The disturbance $w(t) = b(t) + r(t) $ contains the signal $b(t)$ caused by the bumps on the ground and the signal $r(t)$ due to the curvature of the trail. We assume that ratio of the size of $b(t)$ to the size of $r(t)$ is some $\epsilon >0$. The control action $u(t) =  u_\ell(t ) + u_h(t) $ is generated by $u_\ell(t )$ from the reflex loop and $u_h(t)$ from the planning loop. The reflex loop (denote by $L$) compensates for bumps using reflex at a lower-layer, and the planning loop (denoted by $H$) tracks the trail at a higher-layer.

There are speed and accuracy constraints in each control loop. We assume that the reflex loop can transmit signals from sensory to motor units with a signaling rate $R_\ell$ and total delay $T_{\ell} + T_i$, where $T_{\ell}$ models the signaling delay, $T_i$ aggregates other internal delays, and $R_\ell$ and $T_{\ell}$ are subject to a component SAT $R_\ell = \lambda_\ell T_{\ell}$.  The planning loop has a signaling rate $R_h$ and total delay $T_h - T_a$, where $T_h$ models the signaling delay, $T_a$ is the advanced warning, and $R_h$ and $T_{h}$ are subject to a component SAT $R_h = \lambda_h T_{h}$. The difference in their level of internal delay or advanced warning comes from the fact that the control response to trail curvature can be planned in advance by viewing the trail ahead, whereas the bumps are often controlled after a cyclist senses its impact.

With sufficiently large advanced warning $T_a  ( > T_h)$, the state-deviation $ \sup_{ \|  b \|_\infty \leq \epsilon, \|  r \|_\infty \leq 1} \|  x\|_\infty$ is lower-bounded by 
\begin{align}
\label{eq:optimal1}
\left\{  T_\ell +T_i + \frac{  1}{2^{R_\ell} - 1}  \right\}\epsilon
+ \frac{1}{2^{R_h}  -1 }. 
\end{align}
Here, the rates $(R_\ell, R_h)$ are the information capacity used by the subtasks in individual layers, but do not include the information capacity used for other tasks or homeostasis. {\ed This lower bound is tight in the sense that a controller exists that achieves this bound. Analogous to the case of Eq.~\ref{eq:performance_det}, the performance bound in Eq.~\ref{eq:optimal1} holds regardless of whether the planning and reflex layers have feedback or feedforward structures.}

Note that the overall lower-bound for the error is the sum of the errors in the lower reflex layer and the higher visual layer. This property is consistent with experimental observations in Fig.~\ref{fig:experiment-bumps-trail}. This decomposition holds when the bumps $b$ and trail changes $r$ are independent and small enough to be independently controlled by each layer (see Discussion for the situations when this assumption does not hold). Under these assumptions, the feedback control system can be decomposed into two independent subsystems that individually control the deleterious effects of $b$ and $r$. One uses the feedback loop $L$ to control the error dynamics Eq.~\ref{eq:plant} with $r(t) \equiv 0$, while the other uses the feedback loop $H$ to control the error dynamics Eq.~\ref{eq:plant} with $b(t) \equiv 0$. 

The separation of Eq.~\ref{eq:optimal1} into the individual errors caused by two subsystems allows us to use the preceding insight to study the layered control architecture used in the biking tasks. The reflex feedback typically operates in the regime of delayed reaction, as reflexes often sense bumps only after the bike has hit them. The planning feedback typically works in the regime of advanced planning, since the trajectory of the bike and trail can often be seen in advance. From Fig.~\ref{fig:loop-hete-joint-a}, the reflex feedback has the best performance with small signaling delay at the expense of a low signaling rate. On the contrary, the planning feedback has the best performance with a high signaling rate at the expense of a large signaling delay. This theoretical prediction on the relative delays of the two layers parallels the relative delays in bump only and trail only tasks observed in our experiment (see the supplementary material section 4.B) and comply with the existing literature \cite{roland2006cortical,matthews1991human}. 

\begin{figure}[h!]
\centering

\begin{subfigure}[b]{\linewidth}
	\caption{System SATs for visual object tracking}
	\centering
\includegraphics[width=0.8\linewidth]{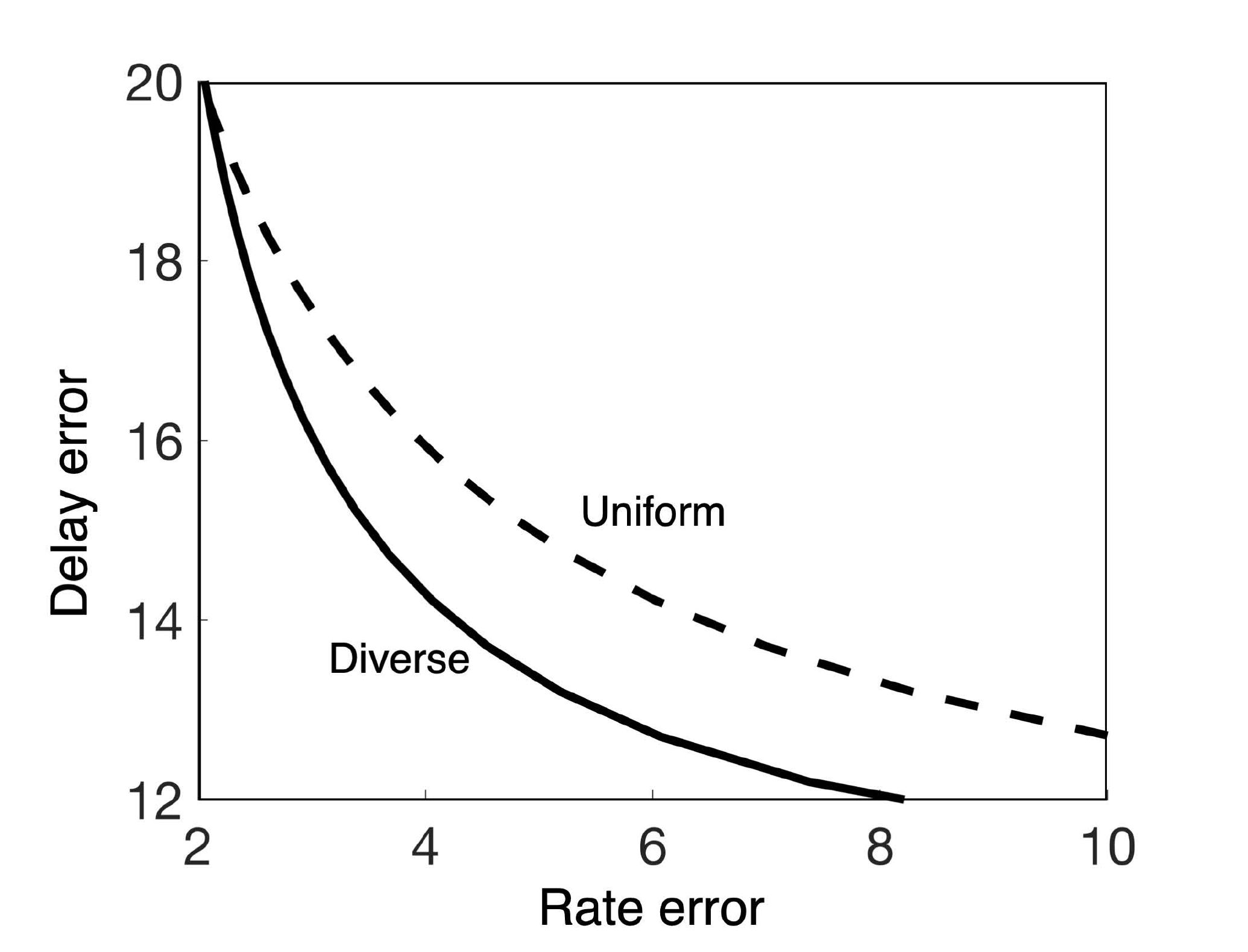} 
\label{fig:DESS-visual}
\end{subfigure}

\begin{subfigure}[b]{\linewidth}
	\caption{System SATs for mountain biking}
	\centering
\includegraphics[width=0.8\linewidth]{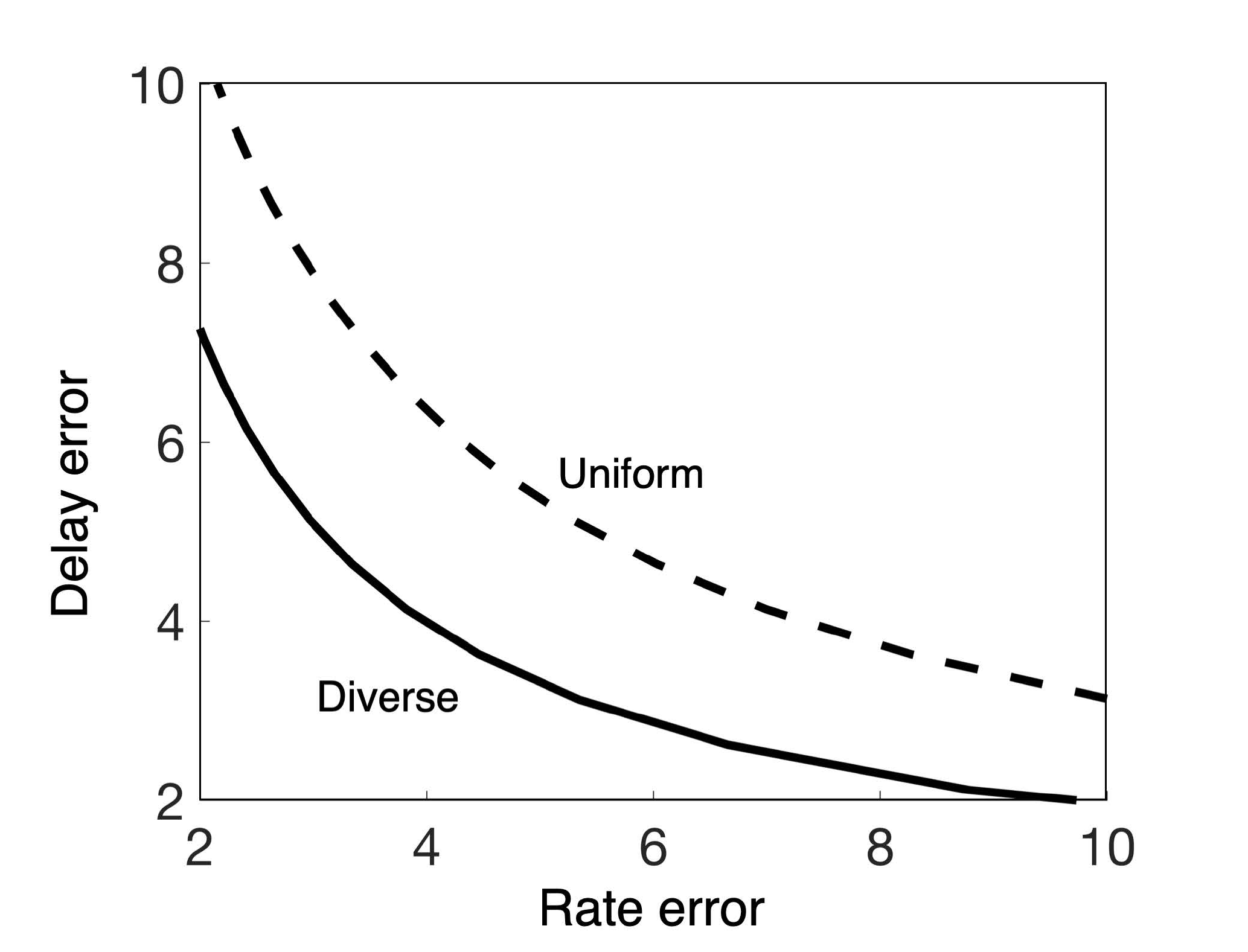}
\label{fig:DESS-bike}
\end{subfigure}

\caption{{\re
\textbf{Diversity in the components (Diverse) improves performance compared with uniform components (Uniform).} In the diverse cases, both layers are allowed to use heterogeneous signaling delays and rates. In uniform cases, they are constrained to be homogeneous. The horizontal axis shows the sum of the rate errors in both layers, and the vertical axis shows the sum of the delay errors in both layers. For both the setting of visual object tracking and lateral control in trail tracking, we can observe Diversity-Enabled Sweet Spots (\ie the diverse cases have less stringent SATs in control than the uniform cases).  
\textit{(A)} For visual object tracking, we used the component SAT $R_{\ell}= 0.1 T_{\ell}$ for the reflex loop, and $R_h = 0.1 T_h$ for the planning loop. The component SATs are converted into system SATs by Eq.~\ref{eq:optimal1} with parameters $T_i=10$, $T_a=10$, and $\epsilon = 1$. Although the plot is shown for specific levels of net delay and net warning, Fig.~\ref{fig:SAT-delay-vs-advanced} suggest that, in the diverse case, the advantageous performance holds over a broad range of net delay/warning as the optimal signaling delay and rate takes a constant value when the net delay and net warning are in $[-4, 10]$ and $[-10, 4]$, respectively. 
\textit{(B)} For the lateral control in trail tracking, we used the component SAT $R_{\ell}= 0.1 T_{\ell}$ for the reflex loop, and $R_h = 0.1 T_h$ for the visual loop. The component SATs are converted into system SATs by Eq.~\ref{eq:optimal1} with parameters $T_i=0$, $T_a=100$, and $\epsilon = 1$. In the diverse case, the reflex layer and the planning layer are allowed to use heterogeneous signaling delays and rates, whereas in the uniform case, they are constrained to be homogeneous.}
} 
\label{fig:DESS}
\end{figure}

The resulting benefit of diversity in delays and rates is illustrated in Fig.~\ref{fig:SAT-system}, which shows the optimal component composition (Fig.~\ref{fig:loop-hete-joint-a}) and compares the system performances of the uniform and diverse cases (Fig.~\ref{fig:loop-hete-joint-b}) when component SATs in Eq~\ref{eq:delay-rate-spike1} are applied into the system performance in  Eq.~\ref{eq:optimal1}. The relaxed system SATs in the diverse case compared to the uniform case suggests that diversity in the layered control architecture helps improve the fundamental performance limits arising from component SATs (Fig.~\ref{fig:DESS-bike}). {\ed When the appropriate diversity and layers do not exist, the disturbance is processed by a control loop whose delay and rate are not optimized for the specifics of its access and the extent of advanced warning. Thus, the uniform case is expected to have a larger error and worse performance, as indicated by the limits in either of the two terms in Eq~\ref{eq:optimal1}, than by control loops with optimized delays and rates.} The diversity between the two layers virtualizes the performance, allowing the overall system to exploit both predictive control and fast reflex to reduce errors.



\section*{Discussion}

Our theoretical analysis of oculomotor control and biking showed that the deleterious effects of component delays and inaccuracies on control performance can be mitigated by layering and diversity. Diversity allows optimal trade-offs between delay error and rate error (Fig.~\ref{fig:DESS}). 


\subsection*{Comparisons with previous studies} At the component level, reducing the energy needed for information transmission is often a major concern~\cite{hasenstaub2010}. At the system level, fast, accurate, and robust control is important for survival~\cite{standage2015toward,franklin2011computational}. Here, we consider the design objective of optimizing the robustness of sensorimotor control given limited biological resources in energy and space. Energy efficiency considers the signaling rate for a component as the design goal -- maximizing information rate given a fixed energy budget -- but in our framework, information transmission is only a means to the goal of efficient control.

\textit{Optimal nerve composition from a system perspective.} The difference in the two design goals leads to different conclusions for the optimal composition of nerves. From the component perspective of maximizing information rate within the energy and resource budget, having many small axons that send information at the lowest acceptable rate is desirable~\cite{laughlin2003communication}. From the system perspective of achieving robust sensorimotor behaviors, balancing speed and accuracy in neural signaling is more important since this minimizes the total control error due to delays and limited rates (Fig.~\ref{fig:SAT}). Conversely, maximizing signaling rate may lead to large delays due to the component SATs, which in turn degrade the robustness in sensorimotor control. These contrasting results reveal the fundamental difference between optimizing component properties and optimizing system performance.

\textit{Enabling factors of robust performance.} Unlike models of sensorimotor control that abstract away the component speed and accuracy limits and assume they are negligible, our model explicitly incorporates these constraints and their impact on system performance using robust networked control theory. Although both types of models can explain the empirical observation that component limits have minimal impact on sensorimotor performance, the factors that enable robust performance arise differently in the two models. 

When the component limits are negligible, to achieve robust performance only requires properly dealing with uncertainties using the mechanisms such as sensorimotor adaptation, optimal feedback control, impedance control, predictive control, Bayesian decision-making, and robust (risk-averse) control~\cite{franklin2011computational,sanger2010controlling}. However, when component limits are not negligible, robust sensorimotor control additionally requires the mechanisms to successfully mitigate the component hard limits through the use of effective layered control architectures with appropriate diversity. In other words, even with a collection of inexpensive layers, which may be slow or inaccurate, an effective layering can virtualize a fast and accurate control response. Such performance is achieved when the layers have proper diversity to collectively span the heterogeneous requirements needed for robustly performing a task  (Fig.~\ref{fig:DESS}). 

\subsection*{Assumptions and limitations} {\re We approached mountain biking as a complex biological system by isolating each part of the complexity in carefully designed experiments and drawing general conclusions that can then be further tested. Our focus was on the diversity of components and layers. We chose an experimental paradigm designed to separate this mechanism from many others that might otherwise have obscured it.}

\textit{Sensorimotor integration.} In the biking task, the subject is stable in the chair, the bumps do not affect visual trail tracking, and the only motor output is the position of the wheel. In this simplified version of biking on a smooth road in an ideal environment, we showed that the errors due to trails and bumps are additive and that the closed-loop performance bound matched the empirical observations. Our basic model captures the essence of how component limitations are deconstrained at the system level and is consistent with the results of our experiments and the properties of nerves. Our model can be extended in several theoretical and experimental directions to better understand the control processes and encompass a broader range of tasks.

\textit{System dynamics.} How do dynamical system properties such as poles and zeros change the relationship between component SATs and the stability and performance of the whole system? The impact of delay in unstable systems may lead to oscillations, a qualitatively different regime from (marginally) stable systems. Such models can be used to model human stick balancing, whose dynamics have unstable poles and zeros.  Pole balancing gets harder when the visual focus shifts towards the lower end of the stick. 

\textit{Disturbances.} What is the impact of specific types of disturbances and their properties on system design? For example, the low rotational inertia and low probability of perturbation of the orbit in oculomotor dynamics may explain why oculomotor control can achieve remarkable performance even with minimal proprioception, which is essential in guiding limb movements and achieving fast and accurate feedback control.

\textit{Sensorimotor control.} 
What other influence do quantization, delay, saturation and other properties have on the performance of motor systems and sensorimotor integration? Component constraints can be further refined by including the specific roles they have in neural coding and their functions in feedback and feedforward pathways~\cite{zhaoping2014understanding,crapse2008corollary,borst1999information,laughlin2003communication,Sepulchre2019}. The model can be extended to account for the constraints of the motor system (muscle strength, accuracy, speed, and fatigue) in the context of reaching~\cite{Nakahira2020Fitts}, throwing~\cite{hore2011skilled}, and biking (\eg a strong biker is often able to recover from a larger displacement than a weaker one). The impact of component constraints on performance in the biking task can be studied for patients with motor disabilities, such as those with Parkinson's disease, who may have disrupted speed and accuracy constraints in their control loops.

\textit{Interactions between layers and control systems.} When biking on a twisty, bumpy road, stabilizing against bumps and visual object tracking is more complex. The VOR and internal feedback loops (particularly within the visual cortex) work together to stabilize vision despite bumps. In this situation, the errors from the two layers may not be additive. Depending on the specifics of the lower layer, the bumps may influence how well the higher layer can sense and make decisions. Including the effects of bike dynamics on head and body movements  could reveal interactions between these controllers and layers.  

In our model, descending prediction errors from sensors were transferred to actuators to take advantage of the available signaling rate. This feedback resembles predictive coding, which may explain the existence of massive feedback from higher cortical areas back to the primary sensory areas of cortex~\cite{muckli2013network,zhaoping2014understanding}.

The further integration of different control loops by matching controller gains, information transmission and coordination between layers is a major challenge. The cerebellum collects proprioceptive and efference copies of motor commands from the entire body and manages gain adaptation while keeping sensorimotor loops from interfering with each other~\cite{coenen1996learning}. For example, in the VOR circuit, vestibular inputs project to both the cerebellar cortex and the vestibular nucleus, and the cerebellum in turn projects back to the vestibular nucleus. The cerebellum adjusts the gain of the feedforward loop using sensory prediction error signals from image slip in the retina during head movements~\cite{shadmehr2010}. Gain is also modulated by verging the eyes. Interestingly, during a vergence eye movement, the gain is reset before the eyes reach the endpoint \cite{Snyder1992, coenen1995Vergence}. For head movements that are not horizontal, the transformation to non-orthogonal eye muscles is even more complex.

\textit{Cognition.} 
More layers that could be added to our model include model-based prediction, memory, cortical representations, alertness and attention, all of which influence computation and communication in the central nervous system. Integrating these additional layers could lead to a better understanding of how distributed control is achieved in brains. Connecting the SATs in sensorimotor control and the SATs in decision making will provide further insights into how control and cognitive processes are optimally integrated.


\textit{Factors that contribute to DESSs.} In our models of oculomotor control and trail following, the higher layer performs predictive planning and control and requires a higher data rate than the lower layer. This is commonly found in engineering systems. For example, a model predictive controller (MPC) or path planner at the higher layer is combined with a PID or robust controller at the lower layer. The optimality condition suggests that the lower layer performs best at a fast time scale, while the higher layer performs best with higher data rates and processing power. Moreover, given an autonomous system that requires heavy computation in decision making, adding the reflex for fast responses can largely reduce the latency requirement of the higher layer. 

Other factors that contribute to DESSs include division of labor~\cite{gyorgy2019brain}, hierarchical analysis of sensory input with different spatial and temporal scales at each successive layer, the diverse properties of the muscles and uncertainties that govern the control of subtasks in each layer.

\textit{Design considerations beyond SATs.} There are other ways to improve system performance. One example is overarm throwing, where moving the arm at high speeds leads to increased accuracy. Although muscle noise increases with force in the normal operating range, it drops at maximum strength. Other examples include the gyroscopic effect: biking transitions from unstable to stable dynamics as the velocity crosses a threshold; figure skate spinning is more stable at a faster speed (gyroscopic effect); hopping helps stabilize balance. This phenomenon could arise from a combination of factors: activating more layers; increased sensing sensitivity and actuation capability due to larger motion amplitude; the use of oscillation to stabilize control; and the gyroscopic effect.

\subsection*{Improving the fundamental performance limits}  DESSs require multiple layers with optimal and robust policies and diversified hardware within each layer.  The overall performance of the layered system can improve on the limits in each layer.  Below, we discuss some of these improvements.

\textit{Multiplexing motor systems and waterbed effects.} An appropriate diversity of layers is needed so that essential tasks can be efficiently multiplexed and performed quickly and accurately. There are many feedback loops in sensorimotor systems that can be multiplexed~\cite{koch2007attention}. We have examined how humans can multitask bike trail tracking and stabilizing against bumps (Fig.~\ref{fig:model biking-system}, Fig.~\ref{fig:experiment-bumps-trail}), but trail tracking while texting --- not anticipated by evolution --- leads to catastrophic crashes.  Robust behavior is accomplished by a deliberate design that has separate layers of sensorimotor pathways for each subtask.

Multiplexing is a costly investment since building and maintaining each layer requires additional biological resources. Moreover, there are waterbed-like effects in the sensitivity of a system to disturbances: In feedback control, suppression in one frequency range necessarily increases disturbances in some other frequency range. This phenomena is captured in Bode's sensivitity integral~\cite{fang2017towards,Leong2017waterbed}. Revisiting Bode's sensitivity integral from the perspective of layering and diversity provides a complementary perspective of disturbance rejection~\cite{fang2017towards,Leong2017waterbed}. 

More generally, waterbed-like effects can occur when designs optimized for one type of environment induce fragility in other types. Thus, in the evolution of layered control architectures with a fixed resource budget, improving the capability for one task may induce fragility in others. Identifying these hidden trade-offs and waterbed-like effects could provide a new evolutionary perspective on the organization of sensorimotor control in brains, effects that have been largely overlooked in most studies. Our framework for analyzing multilayer systems could provide explanations for why some tasks share and others compete for resources.


\textit{Cross-layer learning and optimizing architectures.} DESSs can be achieved by decomposing a new task into subtasks that are implemented in different layers. Each layer has hard limits in speed, accuracy and flexibility in learning and control. Higher layers are often flexible but slow, whereas lower layers are faster but less flexible. Repetitive practice identifies and accumulates evidence on potential subtasks that do not require much flexibility and can be automated in lower layers with improved speed and accuracy. Learning how to efficiently allocate layers allows the system to better virtualize tasks to achieve fast, accurate and flexible behaviors despite layers that are by themselves slow, inaccurate, or rigid.  Suboptimal allocations of layers can expose the hard limits of individual layers to performance bottlenecks. 

This is illustrated in overarm throwing. Beginners often use the central nervous system to think about controlling the release of the ball, but highly-skilled players use feedforward control of finger force and stiffness to control the timing of release~\cite{hore2011skilled}. Feedforward control is faster, more accurate and can rapidly adapt. This is possible because feedforward control operates at a millisecond-level through the exceptional sensitivity of force sensing, which greatly improves throwing accuracy.  Force control has a much lower-dimensional design space for finger muscle stiffness. Low-dimensional rigid motion allows for fast adaptation to new targets and wind conditions. By shifting  high-level control to a lower layer, higher-level resources can be redeployed for other tasks.


Cross-layer learning is qualitatively different from incrementally improving control parameters within a layer. This type of learning happens on a much slower time-scale and leads to abrupt changes in behavior. Because of these properties, cross-layer learning is more difficult to observe and study in controlled experiments. Our layered architecture could serve as a starting point for developing a theory for cross-layer learning~\cite{sejnowski2020unreasonable}.

\subsection*{DESS is a universal design principle}

Diversity is `the most ubiquitous rule' in living systems'~\cite{gyorgy2019brain}. In this paper, we studied the underlying mechanisms through which diversity in the delays and rates of sensing and signaling \textit{between} layers improves control performance. Our companion paper shows that diversity of components \textit{within} a layer also boosts performance, and in particular that Fitts' law for SATs in reaching can be explained by DESSs in motor nerves~\cite{Nakahira2020Fitts}. 

DESSs can also be found in many other systems: Human combines fast and slow decision-making processes~\cite{kahneman2011thinking,trimmer2008mammalian};  The immune system combines fast general responses with slower targeted responses~\cite{smith2018host}; The smart grid combines power flow in a slow layer and frequency control in a fast layer; Internet of Things (IoT) integrates cloud computing (which has high computing capability and are centralized) with edge computing (which can quickly respond to local disturbances).

\section*{Conclusions}

Our case studies are just the tip of the diversity iceberg through which diverse mechanisms in prediction, estimation, and actuation within and between layers boost system performance. Understanding the design principles of layered architectures in biological systems, particularly those that achieve DESSs in delayed, quantized, distributed, and localized control, can inspire the design of robust technological systems, which increasingly face challenges similar to those encountered in human sensorimotor control. Design tools and engineering case studies, in turn, will help distill the design principles in biological systems that enable robust and flexible behaviors through complex and heterogeneous neural mechanisms.



\section*{Materials and Methods} 

We developed a platform for biking games that simulate some aspects of riding a mountain bike~\cite{2019Experimental}.  The platform is inexpensive and easy to implement. During the experiment, the subject looked at a PC monitor and turned a wheel to follow the desired trajectory. The trajectory had a constant velocity for each segment but abruptly switched between right and left segments. The console for the biking task is shown in Fig.~\ref{fig:interface}. We conducted experiments with four participants and recorded their biking trajectories and lateral errors in control. 

To study how layers multiplex, we compared the  behaviors when there are bumps in the road, curvature in the trail and both. In the first task, the bumps were generated by pushing the steering wheel at a constant torque for $0.5$ second. In the second task, the trail was generated with the angle $\theta \in \{10^{\circ}, 20^{\circ},\dots, 80^{\circ} \}$ and alternated between left and right with exponentially distributed time intervals, so that the participants cannot anticipate the abrupt shifts without advanced warning in vision. In the last task, the bump and trail changes were generated independently according to the first two settings. A comparison of the error dynamics of the three tasks is shown in  Fig.~\ref{fig:experiment-bumps-trail}. 

To test of the impact of component SATs, we compared the  behaviors when steering wheel input acts on the position with delays, with quantizations, and both. The worst case errors where measured in the three tasks. In the first task, the added delays were set to be $T=-0.8, -0.6,\cdots,0.2, 0.4$ seconds, where negative delays were realized by adding an advanced warning in the visual input, and the positive delays were implemented by adding an external delay in actuation. In the second task, the rate of the quantizer was set to be $R=1,2,\cdots, 7$ bits per unit time. In the third task, the delay and quantization were added according to the first and second settings, respectively. In addition, the delay and rate are set to satisfy $T=(R-5)/20$, which simulates the component SAT in Eq.~\ref{eq:delay-rate-spike1}. Each set of parameters lasted for $30$ seconds before switching to a new set of parameters. The first $10$ seconds of each $30$ second trial were not used to measure the performance in order to eliminate switching and learning effects. Before each experiment, subjects were trained until their performance stabilized. The errors between the desired and actual trajectory are shown in Fig.~\ref{fig:experiment-SAT}. This plot suggests that the error caused by the added delay and quantization is the sum of the error caused by added quantization and the error caused by added delay, as suggested by the theoretical prediction Eq.~\ref{eq:performance_det} in the deterministic setting. We also tested the average errors in an average-case framework (see Section 4 in the supplementary material). 

\textit{Participants}. All participants gave informed consent. The study protocol number 19-0912 was approved by the Institutional Review Board at the California Institute of Technology.

\textit{Data availability.} All data and programs used to analyze the data are available at https://cnl.salk.edu/\textasciitilde{terry}/DESS-PNAS/.


\acknow {This research was supported by National Science Foundation (NCS-FO 1735004 and 1735003) and the Swartz Foundation. Q.L. was supported by a Boswell fellowship. This paper is based on the theoretical doctoral research of Y.N. and on the experimental research of Q.L.}

\showacknow{} 


\bibliography{PNAS}
\end{document}